# DECOMPOSITION TABLES FOR EXPERIMENTS I. A CHAIN OF RANDOMIZATIONS

By C. J. Brien and R. A. Bailey

*University of South Australia and Queen Mary, University of London*

One aspect of evaluating the design for an experiment is the discovery of the relationships between subspaces of the data space. Initially we establish the notation and methods for evaluating an experiment with a single randomization. Starting with two structures, or orthogonal decompositions of the data space, we describe how to combine them to form the overall decomposition for a single-randomization experiment that is "structure balanced." The relationships between the two structures are characterized using efficiency factors. The decomposition is encapsulated in a decomposition table. Then, for experiments that involve multiple randomizations forming a chain, we take several structures that pairwise are structure balanced and combine them to establish the form of the orthogonal decomposition for the experiment. In particular, it is proven that the properties of the design for such an experiment are derived in a straightforward manner from those of the individual designs. We show how to formulate an extended decomposition table giving the sources of variation, their relationships and their degrees of freedom, so that competing designs can be evaluated.

**1. Introduction.** This paper investigates methods for evaluating designs for experiments with multiple randomizations [14] that follow each other in a chain, like those in Figures 1 and 2. In general, the form of the analysis-of-variance table for an experiment is informative when evaluating a design for it. The foundation of this table is the orthogonal decomposition of the data space that it reflects. Hence the purpose of this paper is to establish the appropriate orthogonal decomposition of the data space for some experiments with multiple randomizations. This requires an extension of the decomposition for standard textbook designs, almost all of which involve a single randomization.

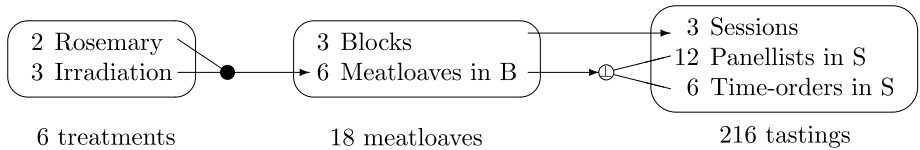

Fig. 1. *Randomization diagram for Example 1: treatments are randomized to meatloaves, which are in turn randomized to tastings;* B *denotes Blocks,* S *denotes Sessions.*

EXAMPLE 1 (Meatloaves). At the Joint Statistical Meetings in New York in 2003, T. B. Bailey described the two-phase sensory experiment shown in Figure 1. This figure is repeated from [14], where the conventions used in such diagrams are explained. In the first phase of the experiment, there were six treatments, consisting of all combinations of a two-level factor Rosemary (present or not) with three quantities of irradiation. These treatments were applied to the process of making meatloaves, using a randomized complete-block design with three blocks. In the second phase, blocks were randomized to sessions. In each session, twelve panellists tasted a portion of all six meatloaves from the block assigned to that session. The design in this phase consisted of a pair of $6 \times 6$ Latin squares in each session. The combined effect of these two randomizations is to randomize treatments to tastings.

The design for each phase of this experiment is orthogonal, in a sense that we make precise in Section 4. It should therefore be straightforward to obtain the appropriate decomposition. However, as T. B. Bailey reported, the analysis of variance given for similar two-phase experiments is frequently wrong.

To obtain the decomposition for experiments with a chain of randomizations, we have to assume that the design for each single randomization satisfies some conditions, which are explained in Section 4. Although they are rather technical, they are satisfied by many designs used in practice. The decomposition that we obtain works just as well for complicated nonorthogonal designs like the one in Figure 2 as it does for orthogonal designs.

Our approach extends that of Fisher [20], Wilk and Kempthorne [38], Nelder [28, 29] and Bailey [1, 9] for a single randomization. The objective of a randomization is to assign one set of objects, here designated $\Upsilon$, to another set of objects, here designated $\Omega$. Frequently $\Upsilon$ is the set of treatments and $\Omega$ is the set of observational units. In Example 1, the set described by the left-hand panel consists of six treatments while the set described by the middle panel consists of 18 meatloaves. The randomization of the 18 meatloaves to the 216 tastings is also considered to be a single randomization, even though there are multiple arrows between their panels. This is because it can be achieved by selecting a single permutation of the 216 tastings, as explained by Brien and Bailey [14].

In Section 3, we define *structure* on a set to be an orthogonal decomposition of the relevant vector space. The structure on the set of treatments in Example 1 consists of the spaces for the grand mean, two main effects and their interaction; the structure on the meatloaves consists of the spaces for the grand mean, differences between blocks and within-block differences. In general, the structure on $\Omega$ is the unrandomized structure, which reflects the topographical, managerial or physical features that are there *before* treatments are assigned.

Although our main results apply to quite general structures, we limit our examples to structures which are uniquely defined by a list of factors, their numbers of levels and their nesting relationships. Such a structure can be succinctly shown in a panel. In the right-hand panel in Figure 1, Panellists and Time-Orders are both nested in Sessions, in the sense of [28] and [14]: thus there are 36 panellists in total, not twelve.

Brien [11] introduced the term *tier* to describe the collection of factors on a set before it is associated with another set by randomization. Most standard textbook designs have two tiers. Brien and Bailey [14] introduced randomization diagrams to summarize such designs. There are two panels, one for each tier. Underneath each tier is written the corresponding set; this may be either a symbol, such as $\Upsilon$ or $\Omega$, or a defining phrase, like "6 treatments." Arrows go from the randomized tier (on $\Upsilon$) to the unrandomized tier (on $\Omega$).

The decomposition for such a design is obtained by starting with the structure on $\Omega$ and then refining it by the structure on $\Upsilon$, according to the relationship between the two structures. The *decomposition* consists of the set of orthogonal subspaces of the data space which is characterized by the set of mutually orthogonally idempotents derived from those of the two structures. It is displayed in a *decomposition table* which has rows and columns. Each row corresponds to one of the subspaces in the decomposition: some authors call this a "source of variation." The number of degrees of freedom shown in this row is just the dimension of the subspace. These tables also show efficiency factors, that we define in Section 4. Bailey [9] gives such decompositions for orthogonal designs, calling the ensuing table the *skeleton analysis of variance*. Most other textbooks on design and analysis of experiments do not do this. Instead, it is usual to write down a model of terms thought to be relevant, which are then listed in an ANOVA table. This method has several disadvantages: the resulting model is somewhat arbitrary [38]; the analysis does not reflect the confounding arising from the randomization; and, for nonorthogonal designs, non-unique partial analyses may be produced [15].

Experiments with two or more randomizations involve three or more tiers, so Brien [11] termed them *multitiered* experiments. Three or more panels are required in their randomization diagrams. As noted in [14], multitiered experiments include two-phase, some superimposed and some single-stage experiments, and some multistage experiments which use the same units at each stage; they do not represent a collection of new designs, but are a class of designs made up of several existing design types. However, general methods for assessing their properties have not been established.

This paper considers those multitiered experiments where the arrows in the randomization diagram follow each other in a chain. This commonly occurs in two-phase experiments [27], particularly those that include a second, laboratory phase after an initial phase: examples are sensory experiments, field experiments followed by laboratory processing, and gene-expression microarray studies. In the simplest case, there are two randomizations and three sets of objects, $\Omega$, $\Upsilon$ and $\Gamma$, each with a tier of factors and a structure: $\Gamma$ is randomized to $\Upsilon$, and $\Upsilon$ is randomized to $\Omega$. Each randomization may employ a standard textbook design. In [14], such a pair of randomizations is called *randomized-inclusive* if the systematic design for the allocation of $\Upsilon$ to $\Omega$ must use information from the outcome of the randomization of $\Gamma$ to $\Upsilon$; otherwise the pair is *composed*, as in Example 1.

Section 2 uses a nonorthogonal two-tiered experiment to demonstrate the approach of using the decomposition table to evaluate the design. It then poses the problem of how to obtain the decomposition table for a three-tiered experiment. Section 3 makes precise what we mean by structure on a single set of objects. Section 4 outlines the relevant concepts for an experiment with a single randomization, in order to make them precise and to establish notation. In particular, we show how to take two structures and combine them to produce a single decomposition table for a two-tiered experiment. Sections 5 and 6 extend these results to obtain the decomposition table for experiments involving a pair of randomizations that are either composed or randomized-inclusive. Section 7 generalizes this to longer chains of randomizations.

The analysis-of-variance table for the analysis of a response variable is an extension of the decomposition table. It includes a sum of squares in each row. This has the form $\mathbf{y'Dy}$ where $\mathbf{y}$ is the data vector and $\mathbf{D}$ is the idempotent matrix of rank $d$ which projects onto the subspace in the decomposition for the row. The subspace is variously called the *image* of $\mathbf{D}$, written $\operatorname{Im}\mathbf{D}$, or the *column space* of $\mathbf{D}$. The mean square, $\mathbf{y'Dy}/d$, is usually shown. Sometimes the expected value of this mean square, based on either a randomization model or a mixed model for the response variable, is also included. In this paper, we are not attempting to obtain the whole analysis-of-variance table, merely the appropriate orthogonal decomposition of the data space so that the properties of an experimental design can be established. In particular, we have not given the expected mean squares, which are often of use in assessing a design. However, for the randomization model, they are easily deduced from the decomposition table, and so the likely relative sizes of variances of effects under this model can be gauged.

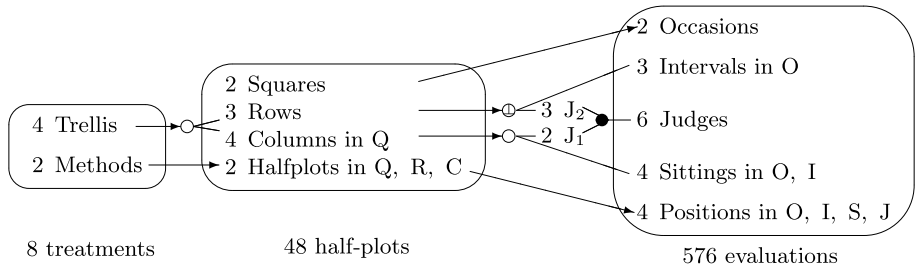

FIG. 2. *Randomization diagram for Example 2: treatments are randomized to half-plots, which are, in turn, randomized to evaluations; Q, R, C, O, I, S, J denote Squares, Rows, Columns, Occasions, Intervals, Sittings and Judges, respectively.*

**2. Evaluating the design for one experiment.** In designing a standard two-tiered experiment, one considers the sources of natural variation, the sources with which the treatments are confounded and the amount of information about different treatment contrasts available from the different sources.

EXAMPLE 2 (A viticultural experiment). Brien and Payne [15] describe a viticultural experiment comprising the field phase of a two-phase experiment. Two adjacent Youden squares are used to assign trellis treatments to the plots, a plot being a row-column combination within a square. Each plot is divided into two half-plots, to which two methods of pruning are assigned at random. There are two sets of objects involved in this phase, half-plots and treatments, with corresponding panels shown in the middle and left of Figure 2. In this example, Q will be the single letter used for Squares to avoid a conflict with Sittings.

When this phase of the experiment is being planned, the next step is to establish properties of the design by obtaining the decomposition in Table 1 (further details about its derivation are in Section 4). The left-hand side of this table contains what Fisher [20] called the "topographical analysis;" it reflects the unrandomized structure on the half-plots. The other half of the table gives the randomized structure derived from the treatments. The table also shows the relationship between the two structures. Trellis is partially confounded with the sources Columns[Squares] and Rows#Columns[Squares], with efficiency factors 1/9 and 8/9, respectively. Method and Trellis#Method are both totally confounded with Halfplots[Squares ∧ Rows ∧ Columns]. Hence, in a randomization model, the expected mean square for the Residual in Rows#Columns[Squares], with nine degrees of freedom, is equal to the variance $\xi_{\text{RC}[Q]}$ of the Row–Column interaction within Squares. Further, the variance of a difference between two types of trellis, if estimated solely from contrasts in Rows#Columns[Squares], that is, only from contrasts between Row–Column combinations within Squares that are orthogonal to contrasts for both Rows and Columns within Squares, is $(2/12) \times (9/8)\xi_{\text{RC}[Q]}$. Then the designer may be satisfied that most of the information about types of trellis is confounded with a source from which many larger sources of variation have been eliminated.

In the subsequent sensory phase, the half-plots from the first phase are randomized, using two Latin squares and an extended Youden design, to glasses in positions on a table for evaluation by judges. This sensory phase adds a third set of objects, evaluations, with associated panel shown on the right of Figure 2. What are the properties of this two-phase experiment, whose randomizations are composed? We would like to establish the decomposition table for the full experiment. How should one combine the properties of the second-phase design with those of the first-phase design to yield the decomposition table for the final two-phase design?

TABLE 1
*Decomposition table for the first phase of Example 2*

| half-plots tier | | | treatments tier | |
|---|---|---|---|---|
| source | d.f. | eff. | source | d.f. |
| Mean | 1 | 1 | Mean | 1 |
| Rows | 2 | | | |
| Squares | 1 | | | |
| Rows#Squares | 2 | | | |
| Columns[Squares] | 6 | $\frac{1}{9}$ | Trellis | 3 |
| | | | Residual | 3 |
| Rows#Columns[Squares] | 12 | $\frac{8}{9}$ | Trellis | 3 |
| | | | Residual | 9 |
| Halfplots[Squares ∧ Rows ∧ Columns] | 24 | 1 | Method | 1 |
| | | 1 | Trellis#Method | 3 |
| | | | Residual | 20 |

Decomposition tables conveniently summarize the properties of designs, so they are useful in evaluating potential designs. Such tables go back to, at least, Fisher [21], Table 29. They are closely related to the *skeleton analysis-of-variance* tables that are given in [9] and that GenStat [32] produces in response to the ANOVA directive with no data. In GenStat tables, sources are relabelled to include only original factors (not pseudofactors), along with their interacting and nesting factors, while the efficiency factors (other than 0 and 1) are listed below the table. Similar tables can be produced with the programs R [36] and S-Plus [25] by analysing randomly generated data. These analysis-of-variance tables have separate subtables for each unrandomized source, such as sources for the half-plots tier in Example 2; these sources are labelled *strata* in GenStat and [9]. Federer [19] also gives separate subtables, with headings such as "main-plot analysis."

The pair of randomizations in Example 2 form a chain. In what order should we do the decompositions for such randomizations: (a) start with the decomposition of the half-plots tier, refine it by the decomposition of the treatments tier as in Table 1 and then use this refinement to refine the decomposition of the evaluations tier, or (b) start with decomposition of the evaluations tier, refine it by the decomposition of the half-plots tier and then refine this refinement using the decomposition of the treatments tier? These questions are addressed in full in Sections 5 and 6.

**3. Structure in a single tier.** Given a set $\Omega$ of objects, let $V_\Omega$ be the space of all real vectors indexed by $\Omega$. This space is also written as $\mathbb{R}^\Omega$ and its dimension is $|\Omega|$. If the elements of $\Omega$ are $\omega_1, \ldots, \omega_n$, then a vector in $V_\Omega$ has the form $(v_{\omega_1}, \ldots, v_{\omega_n})$. The advantages of using $\Omega$ to index the space $V_\Omega$ are: (i) we do not need to agree on a convention for the order of writing the elements of $\Omega$; and (ii) if two sets $\Omega$ and $\Upsilon$ have the same size, we can still distinguish between $V_\Omega$ and $V_\Upsilon$ because they are indexed by different labels.

Following Bailey [4], we interpret *structure* on $\Omega$ to mean an orthogonal decomposition of the space $V_\Omega$. Such a decomposition is specified by the set $\mathcal{P}$ of symmetric, idempotent, mutually orthogonal matrices projecting onto the subspaces of $V_\Omega$ in the decomposition. Thus if $\mathbf{P} \in \mathcal{P}$ then $\mathbf{P}$ is symmetric and idempotent; if $\mathbf{P}_1$ and $\mathbf{P}_2$ are in $\mathcal{P}$ with $\mathbf{P}_1 \neq \mathbf{P}_2$ then $\mathbf{P}_1 \mathbf{P}_2 = \mathbf{0}$; and $\sum_{\mathbf{P} \in \mathcal{P}} \mathbf{P} = \mathbf{I}_\mathcal{P}$, the identity matrix on $V_\Omega$. We shall identify such a set $\mathcal{P}$ of matrices with the corresponding orthogonal decomposition of $V_\Omega$. Each of these matrices has $|\Omega|$ rows and $|\Omega|$ columns, which are labelled by the elements of $\Omega$, so it is called an $\Omega \times \Omega$ matrix. As explained in [7], this means that we do not have to agree on an order for writing down the elements of $\Omega$, so long as we show the labels for the rows and columns.

In our examples, we restrict attention to structures in which each set is uniquely indexed by a set of factors, referred to as a tier. The notation and conventions for such tiers are given in [14], Section 2.2. We use $F_1 \wedge \cdots \wedge F_n$ to denote the *generalized factor* whose levels are the levels combinations of $F_1, F_2, \ldots$ and $F_n$, for $n \geq 1$. Those with $n > 1$ are called *joint factors* by Hinkelmann and Kempthorne [23], Section 4.12.1. If a generalized factor contains $F_i$ and is to be meaningful, then it must also contain every $F_j$ which nests $F_i$. We consider only such generalized factors, calling them *intrinsic*. In addition, we use the trivial factor Universe, which has a single level.

The factors in a tier are partially ordered by the nesting relationship in [28], so they form a partially ordered set, or poset [10]. When there is one object in the set of objects for each combination of the levels of these factors, then the ordered pair consisting of the set of objects and the collection of intrinsic generalized factors derived from the poset is called a *poset block structure* [2, 3, 6, 7]. In this case, there is one idempotent for each intrinsic generalized factor, and we call this collection of idempotents a *poset structure*.

Sometimes, in addition to the tier, pseudofactors are needed to describe the randomization. (The introduction of pseudofactors has no effect on the size of the set of objects.) As in [14], if $F$ is a factor or generalized factor, then each level of a pseudofactor for $F$ corresponds to one or more levels of $F$. The levels of a pseudofactor have no inherent interest. Pseudofactors can be combined with other factors and pseudofactors into *generalized pseudofactors*. In these cases, $\mathcal{P}$ may include idempotents derived from the pseudofactors.

The generalized factors and pseudofactors can be shown on a Hasse diagram [2, 3, 8, 9, 37], that shows the important marginality relationships between them [6, 12, 31, 33, 37]. For generalized (pseudo)factors $G$ and $H$, we say that $H$ is *marginal* to $G$ if each level of $G$ occurs with only one level of $H$ but $G$ is not aliased with $H$. For poset block structures, $H$ is marginal to $G$ if the factors in $H$ are a proper subset of those in $G$. In the Hasse diagram there is a small circle for each generalized (pseudo)factor $G$; the number $n_G$ of levels of $G$ is shown beside it. The full factor name is sometimes abbreviated to its first letter so long as it has been given in full at a higher point in the diagram. If $H$ is marginal to $G$ then $H$ is drawn higher than $G$ and joined to it by downward lines. At the top of the diagram is the factor Universe, which is marginal to every other generalized factor.

Corresponding to each generalized (pseudo)factor $G$ is an averaging operator $\mathbf{A}_G$. This is an $\Omega \times \Omega$ matrix whose $(\alpha, \beta)$-entry is 0 if the objects $\alpha$ and $\beta$ have different levels of $G$ and otherwise is the reciprocal of the replication of the shared level of $G$. Thus, if $\mathbf{y} \in V_\Omega$, then the $\omega$-entry of $\mathbf{A}_G \mathbf{y}$ is the average value of $y_\alpha$ for elements $\alpha$ of $\Omega$ with the same level of $G$ as $\omega$. In particular, $\mathbf{A}_{\text{Universe}} = |\Omega|^{-1} \mathbf{J}$ where $\mathbf{J}$ is the $\Omega \times \Omega$ all-ones matrix.

Also associated with each generalized factor $G$ is a *source*: this is the subspace of $V_\Omega$ for differences between the levels of $G$ that are not accounted for by generalized factors marginal to $G$. Similarly, there is also a *pseudosource* associated with each generalized pseudofactor. Many different notations are in use for nesting and interacting factors in sources—Heiberger [22], Section 12.4, compares several. Here we use a notation that is reasonably intuitive and is unambiguous for poset block structures. In such a structure, the factors in a generalized factor $G = F_1 \wedge \cdots \wedge F_n$ can always be ordered in such a way that (i) if $F_j$ is nested in $F_i$ then $j > i$; and (ii) there is some $m < n$ such that $F_i$ nests another factor in $G$ if and only if $i \leq m$. We write the source for $G$ as $F_{m+1}\# \cdots \#F_n[F_1 \wedge \cdots \wedge F_m]$ if $m \geq 1$; otherwise as $F_1\# \cdots \#F_n$.

The collection of sources and pseudosources specify an orthogonal decomposition of $V_\Omega$, with an orthogonal projector $\mathbf{P}$ in $\mathcal{P}$ for each (pseudo)source. Starting with $\mathbf{P}_{\text{Universe}} = \mathbf{A}_{\text{Universe}}$, each $\mathbf{P}_G$ is obtained by subtracting from $\mathbf{A}_G$ all those $\mathbf{P}_H$ for which $H$ is marginal to $G$. Thus the effect of any $\mathbf{P}$ on any vector is achieved by a straightforward sequence of averaging operations, each equivalent to applying one of the $\mathbf{A}$ matrices, and subtractions. Degrees of freedom $d_G$ are calculated from the numbers $n_G$ by exactly the same process; they are also shown in the Hasse diagram, with their corresponding sources.

EXAMPLE 3. If $\Omega$ consists of $bp$ plots grouped into $b$ blocks of equal size, then the panel for it is like the right-hand one in [14], Figure 1, similar to the middle panel of Figure 1. The poset of factors for plots is the tier {Blocks, Plots} with the $p$-level factor Plots nested in the $b$-level factor Blocks. The set of intrinsic generalized factors derived from this poset is {Universe, Blocks, Blocks ∧ Plots}, and the factor Blocks ∧ Plots, whose levels are the combinations of levels of Blocks and Plots, uniquely indexes the plots. Hence the set of plots and this set of generalized factors form a poset block structure. The generalized factors, with their numbers of levels, are given on the left in the left-hand Hasse diagram in Figure 3. Also, on the right in this diagram are the three sources in the orthogonal decomposition of $V_{\text{plots}}$, deduced from the generalized factors, with their degrees of freedom; only the first letter of each factor name is shown. The source Plots[Blocks] is the subspace for Plot differences within Blocks. It is obtained as differences between all combinations of Plots and Blocks after allowing for the subspace Blocks. The poset structure on plots is the orthogonal decomposition of the $bp$-dimensional $V_{\text{plots}}$ specified by the set of idempotents $\mathcal{P} = \{\mathbf{P}_{\text{Mean}}, \mathbf{P}_{\text{Blocks}}, \mathbf{P}_{\text{Plots[B]}}\}$.

The Hasse diagram on the right of Figure 3, formed using the same process as that on the left, gives the expressions for the elements of $\mathcal{P}$. It shows that $\mathbf{P}_{\text{Mean}} = \mathbf{A}_{\text{Universe}}$, $\mathbf{P}_{\text{Blocks}} = \mathbf{A}_{\text{Blocks}} - \mathbf{A}_{\text{Universe}}$, and $\mathbf{P}_{\text{Plots[Blocks]}} = \mathbf{A}_{\text{Blocks} \wedge \text{Plots}} - \mathbf{P}_{\text{Blocks}} - \mathbf{P}_{\text{Mean}} = \mathbf{A}_{\text{Blocks} \wedge \text{Plots}} - \mathbf{A}_{\text{Blocks}} = \mathbf{I} - \mathbf{A}_{\text{Blocks}}$.

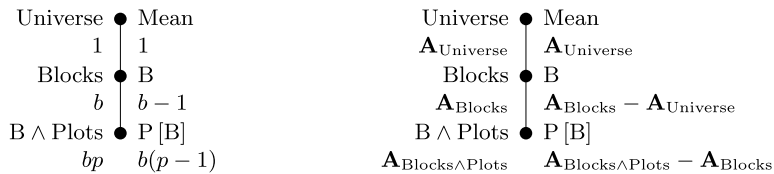

Fig. 3. *Hasse diagram for Example 3: on the left, it is used to calculate degrees of freedom from numbers of levels of generalized factors; on the right, it is used to calculate idempotents from averaging matrices.*

Suppose that $b = 2$ and $p = 3$. Then $|\Omega| = 6$ and let $\Omega = \{\omega_1, \omega_2, \omega_3, \omega_4, \omega_5, \omega_6\}$ be the set of plot labels. Further suppose that the levels of Blocks $\wedge$ Plots for the plots are 11, 12, 13, 21, 22, 23, respectively, where the first number indicates the Block and the second the Plot. In this case the $\Omega \times \Omega$ Block averaging operator is

$$\mathbf{A}_{\text{Blocks}} = \begin{array}{c} \\ \omega_1 \\ \omega_2 \\ \omega_3 \\ \omega_4 \\ \omega_5 \\ \omega_6 \end{array} \begin{array}{c} \omega_1\ \omega_2\ \omega_3\ \omega_4\ \omega_5\ \omega_6 \end{array} \left[ \begin{array}{cccccc} \frac{1}{3} & \frac{1}{3} & \frac{1}{3} & 0 & 0 & 0 \\ \frac{1}{3} & \frac{1}{3} & \frac{1}{3} & 0 & 0 & 0 \\ \frac{1}{3} & \frac{1}{3} & \frac{1}{3} & 0 & 0 & 0 \\ 0 & 0 & 0 & \frac{1}{3} & \frac{1}{3} & \frac{1}{3} \\ 0 & 0 & 0 & \frac{1}{3} & \frac{1}{3} & \frac{1}{3} \\ 0 & 0 & 0 & \frac{1}{3} & \frac{1}{3} & \frac{1}{3} \end{array} \right].$$

EXAMPLE 4. If the set of eight treatments is indexed by all combinations of four trellis treatments with two methods of pruning, then there is no nesting. The poset of factors in the panel on the left of Figure 2 gives the Hasse diagram in Figure 4. This shows the four intrinsic generalized factors and their numbers of levels that, together with the set of treatments, form a poset block structure. It also gives the four sources derived from the generalized factors, and their degrees of freedom. The source Trellis#Methods is the subspace for the interaction of Trellis and Methods. It is obtained as differences between all combinations of Trellises and Methods after allowing for the subspaces Trellis and Methods. The poset structure on treatments is the orthogonal decomposition of the 8-dimensional vector space $V_{\text{treatments}}$ specified by the set of idempotents $\mathcal{Q} = \{\mathbf{Q}_{\text{Mean}}, \mathbf{Q}_{\text{Trellis}}, \mathbf{Q}_{\text{Method}}, \mathbf{Q}_{\text{T\#M}}\}$.

**4. Structure for experiments with a single randomization.** Here we recall the results and concepts from Nelder [28, 29] and James and Wilkinson [26] for an experiment with a single randomization of a set $\Upsilon$ of treatments onto a set $\Omega$ of observational units; we also establish some notation. Although we speak of randomizing the treatments to the observational units, the allocation of treatments to units is actually a function $f$ from $\Omega$ to $\Upsilon$:

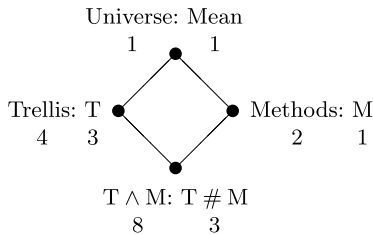

Fig. 4. *Hasse diagram for Example 4, used to calculate degrees of freedom.*

the treatment on observational unit $\omega$ is just $f(\omega)$ [1]. This is shown diagrammatically in Figure 5.

The design function $f$ is the result of deliberate randomization. As discussed in [14], Section 2, the randomization is often achieved by applying a single permutation (perhaps chosen from a restricted set) to the (labels of the) observational units, after an initial systematic plan has been written down. The initial plan is carefully *chosen* by the experimenter from among many possibilities. Part of his or her choice is governed by practical considerations such as the limiting of levels of one treatment factor to large parts of $\Omega$ like whole plots; part is governed by combinatorial considerations such as making an incomplete-block design balanced; part is concerned with getting desirable values for the efficiency factors, which are discussed in the remainder of this paper.

As in Section 3, there is a structure $\mathcal{P}$ on the $|\Omega|$-dimensional vector space $V_\Omega$. Similarly, the structure on the $|\Upsilon|$-dimensional vector space $V_\Upsilon$ may be represented by a set $\mathcal{Q}$ of symmetric, idempotent, mutually orthogonal, $\Upsilon \times \Upsilon$ matrices which sum to $\mathbf{I}_\mathcal{Q}$, the identity matrix on $V_\Upsilon$. The basic problem is to obtain a single decomposition of $V_\Omega$ that combines the two structures $\mathcal{P}$ and $\mathcal{Q}$ to produce a single set of orthogonal idempotents that we denote by $\mathcal{P} \triangleright \mathcal{Q}$.

We call $\mathcal{P} \triangleright \mathcal{Q}$ the set of idempotents for "$\mathcal{P}$ refined by $\mathcal{Q}$" or "$\mathcal{P}$ decomposed by $\mathcal{Q}$." It will be shown that, for structure-balanced experiments, defined below, the elements of this set are of two types: $\mathbf{P} \triangleright \mathbf{Q}$, called "$\mathbf{P}$ pertaining to $\mathbf{Q}$," and $\mathbf{P} \vdash \mathcal{Q}$, called "$\mathbf{P}$ orthogonal to $\mathcal{Q}$" or "the residual of $\mathcal{Q}$ in $\mathbf{P}$." For each $\mathbf{P}$ in $\mathcal{P}$, these decompose $\operatorname{Im} \mathbf{P}$, the subspace of $V_\Omega$ onto which $\mathbf{P}$ projects, into subspaces of the form $\operatorname{Im} \mathbf{PQ}$, for $\mathbf{Q} \in \mathcal{Q}$, and a sub-space that is orthogonal to all of the subspaces $\operatorname{Im} \mathbf{Q}$. If, for some particular $\mathbf{P}$, all of the subspaces $\operatorname{Im} \mathbf{Q}$ are orthogonal to $\operatorname{Im} \mathbf{P}$, then $\mathbf{P} \vdash \mathcal{Q} = \mathbf{P}$.

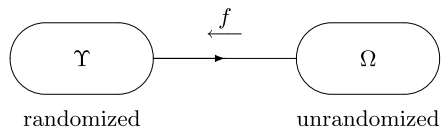

Fig. 5. *Diagram of an experiment with a single randomization.*

However, a little care is needed here because, in general, $|\Omega| \geq |\Upsilon|$, and so the size of a $\mathbf{P} \in \mathcal{P}$ can differ from that of a $\mathbf{Q} \in \mathcal{Q}$. Clearly, we would like all matrices to be $\Omega \times \Omega$ matrices. To achieve this requires the mapping of the space $V_\Upsilon$ into an appropriate subspace of $V_\Omega$. Now, the function $f$ is often represented by the $\Omega \times \Upsilon$ design matrix $\mathbf{X}$, whose $(\omega, t)$ entry is equal to 1 if $f(\omega) = t$ and to 0 otherwise. The column space of $\mathbf{X}$, which we shall temporarily call $V_\Upsilon^f$, is the required subspace of $V_\Omega$; it is isomorphic to the space $V_\Upsilon$.

Next we have to turn $\mathcal{Q}$ into a structure $\mathcal{Q}^f$ on $V_\Upsilon^f$. If generalized factor $G_i$ has averaging matrix $\mathbf{A}_i$ on $V_\Upsilon$, then its averaging matrix $\mathbf{A}_i^f$ on $V_\Omega$ is $\mathbf{X}(\mathbf{A}_i \mathbf{X}' \mathbf{X} \mathbf{A}_i)^- \mathbf{X}'$. In order for generalized factors $G_i$ and $G_j$ on $\Upsilon$ to remain orthogonal when regarded as factors on $\Omega$, we must have $\mathbf{A}_i^f \mathbf{A}_j^f = \mathbf{A}_j^f \mathbf{A}_i^f$; that is,

$$\begin{aligned}(4.1) \quad &\mathbf{X}(\mathbf{A}_i \mathbf{X}' \mathbf{X} \mathbf{A}_i)^- \mathbf{X}' \mathbf{X} (\mathbf{A}_j \mathbf{X}' \mathbf{X} \mathbf{A}_j)^- \mathbf{X}' \\ &= \mathbf{X}(\mathbf{A}_j \mathbf{X}' \mathbf{X} \mathbf{A}_j)^- \mathbf{X}' \mathbf{X} (\mathbf{A}_i \mathbf{X}' \mathbf{X} \mathbf{A}_i)^- \mathbf{X}'.\end{aligned}$$

This is equivalent to the well-known "proportional meeting" condition [6, 9, 37], which is always satisfied if $f$ is equireplicate. We assume that (4.1) holds. Then, like $\mathbf{P}$ in $\mathcal{P}$, the projectors $\mathbf{Q}^f$ in $\mathcal{Q}^f$ can be calculated from the $\mathbf{A}^f$ by suitable subtraction. However, if $f$ is not equireplicate then $\mathbf{Q}^f$ is not, in general, equal to $\mathbf{X}(\mathbf{Q}\mathbf{X}'\mathbf{X}\mathbf{Q})^- \mathbf{X}'$.

The matrices $\mathbf{Q}^f$, for $\mathbf{Q}$ in $\mathcal{Q}$, are symmetric, mutually orthogonal idempotents which sum to $\mathbf{I}_\mathcal{Q}^f$. However, $\mathbf{I}_\mathcal{Q}^f$ is not the $\Omega \times \Omega$ identity matrix: for $\mathbf{w}$ in $V_\Omega$, we have $\mathbf{I}_\mathcal{Q}^f \mathbf{w} = \mathbf{w}$ if $\mathbf{w} \in V_\Upsilon^f$, but $\mathbf{I}_\mathcal{Q}^f \mathbf{w} = \mathbf{0}$ if $\mathbf{w}$ is orthogonal to $V_\Upsilon^f$. In fact, $\mathbf{I}_\mathcal{Q}^f = \mathbf{X}(\mathbf{X}'\mathbf{X})^- \mathbf{X}'$, but we prefer to use notation that suggests that this matrix is an identity for certain other matrices.

Henceforth we identify $V_\Upsilon$ with $V_\Upsilon^f$, $\mathcal{Q}$ with $\mathcal{Q}^f$ and $\mathbf{Q}$ with $\mathbf{Q}^f$, on the understanding that it is the particular choice of $f$ that specifies which subspace of $V_\Omega$ is to be regarded as $V_\Upsilon$. Thus we now regard each $\mathbf{Q}$ as an $\Omega \times \Omega$ matrix. The matrix $\mathbf{I}_\mathcal{Q}$ is a multiplicative identity for each $\mathbf{Q}$, even though it is not for every $\Omega \times \Omega$ matrix.

Recall that $\Omega$ also has a structure $\mathcal{P}$. To be able to combine $\mathcal{P}$ and $\mathcal{Q}$ into a single structure in a unique way, we need $\mathcal{P}$ and $\mathcal{Q}$ to satisfy the condition of structure balance in Definition 1 below. This definition builds on the ideas of [29, 30], in which the design is called *generally balanced* if there are scalars $\lambda_{\mathbf{PQ}}$ for $\mathbf{P}$ in $\mathcal{P}$ and $\mathbf{Q}$ in $\mathcal{Q}$ such that

$$(4.2) \quad \mathbf{QPI}_\mathcal{Q} = \lambda_{\mathbf{PQ}} \mathbf{Q}$$

for all $\mathbf{P}$ in $\mathcal{P}$ and all $\mathbf{Q}$ in $\mathcal{Q}$. (The $\lambda_{\mathbf{PQ}}$ necessarily lie between 0 and 1.) This is equivalent to the conjunction of the following two conditions:

$$\mathbf{QPQ} = \lambda_{\mathbf{PQ}} \mathbf{Q} \tag{4.3}$$

for all $\mathbf{P}$ in $\mathcal{P}$ and all $\mathbf{Q}$ in $\mathcal{Q}$, and

$$\mathbf{Q}_1 \mathbf{P} \mathbf{Q}_2 = \mathbf{0} \tag{4.4}$$

for all $\mathbf{P}$ in $\mathcal{P}$ and all $\mathbf{Q}_1$, $\mathbf{Q}_2$ in $\mathcal{Q}$ with $\mathbf{Q}_1 \neq \mathbf{Q}_2$.

For each fixed $\mathbf{Q}$ and $\mathbf{P}$, condition (4.3) says that every vector $\mathbf{w}$ in the space $\operatorname{Im} \mathbf{Q}$ makes an angle $\cos^{-1} \sqrt{\lambda_{\mathbf{PQ}}}$ with the space $\operatorname{Im} \mathbf{P}$. The $\lambda_{\mathbf{PQ}}$ are called *canonical efficiency factors* in [26]; we abbreviate this to *efficiency factor*. For block and row–column designs, they are the same as those obtained using the reduced normal equations [24, 26].

For each fixed $\mathbf{P}$, condition (4.4) says that the orthogonal spaces $\operatorname{Im} \mathbf{Q}_1$ and $\operatorname{Im} \mathbf{Q}_2$ remain orthogonal when projected onto $\operatorname{Im} \mathbf{P}$. It is a form of adjusted orthogonality [18].

For a single $\mathbf{P}$ and $\mathbf{Q}$, James and Wilkinson [26] said that $\mathbf{Q}$ has *first-order balance* in relation to $\mathbf{P}$ if condition (4.3) holds. However, they did not require that $\mathbf{P}$ and $\mathbf{Q}$ be members of a set of mutually orthogonal, symmetric idempotents.

In [29, 30], general balance is defined for the case that $\mathcal{P}$ and $\mathcal{Q}$ are simple orthogonal block structures, which are those poset block structures that can be defined by formulas using crossing and nesting. In their expository paper, Houtman and Speed [24] took $\mathcal{P}$ as given, although not necessarily defined by factors, and defined the design to be generally balanced if there exists any orthogonal decomposition $\mathcal{Q}$ of $V_\Upsilon$ such that (4.2) holds. This offended writers such as Pearce [35], who do not like the idea that some designs, such as all incomplete-block designs, can be made generally balanced by using a $\mathcal{Q}$ that is a meaningless decomposition. Payne and Tobias [33] expressed a strong preference for general balance to be a property of a specific $\mathcal{Q}$. Brien [13] defined the design $f$ to be *structure balanced* if, for given structures $\mathcal{P}$ and $\mathcal{Q}$, condition (4.2) holds; Bailey [5] called this general balance *with respect to* $\mathcal{P}$ and $\mathcal{Q}$. They had in mind that $\mathcal{P}$ and $\mathcal{Q}$ would often be determined by factor and pseudofactor relations as described in Section 3, which is the case for all examples in this paper.

The randomization of $\Upsilon$ to $\Omega$ includes the random choice of one of the permutations of $\Omega$ which preserve the structure $\mathcal{P}$ on $\Omega$. Thus $\mathcal{P}$ is inherent in the randomization. However, there can be a choice of structure on $\Upsilon$. It may be that the natural structure $\mathcal{Q}_1$, as defined by factors and marginality, is not itself structure balanced in relation to $\mathcal{P}$, but has a refinement $\mathcal{Q}_2$ involving pseudosources that is. If $\mathbf{Q}$ in $\mathcal{Q}_1$ is partly confounded with more than one $\mathbf{P}$, then, either $\mathbf{Q}$ has first-order balance in relation to every $\mathbf{P}$, or $\mathbf{Q}$ needs to be decomposed into pseudosources $\mathbf{Q}_1, \ldots, \mathbf{Q}_m$, say, in such a way that each $\mathbf{Q}_i$ has first-order balance in relation to every $\mathbf{P}$. Pseudosources play the important role of decomposing some sources into pseudosources and a remainder, so that $\mathcal{Q}_2$ is structure balanced in relation to $\mathcal{P}$. Compared to other possible decompositions of the sources, the advantage of pseudosources defined by pseudofactors is that all projections can still be obtained by averaging and subtracting. Some such pseudofactors are used in randomizing the design (as in the second phase of Examples 2 and 5); others are introduced only as an aid to the analysis (as in a group-divisible block design).

For clarity, we shall use the following definitions, which apply to all structures, whether or not they are natural, and whether or not they are defined by factors and pseudofactors.

DEFINITION 1. A structure $\mathcal{Q}$ is *structure balanced* in relation to a structure $\mathcal{P}$ if there are scalars $\lambda_{\mathbf{PQ}}$ for $\mathbf{P}$ in $\mathcal{P}$ and $\mathbf{Q}$ in $\mathcal{Q}$ such that

(i) $\mathbf{QPQ} = \lambda_{\mathbf{PQ}}\mathbf{Q}$ for all $\mathbf{P}$ in $\mathcal{P}$ and all $\mathbf{Q}$ in $\mathcal{Q}$, and
(ii) $\mathbf{Q}_1\mathbf{PQ}_2 = \mathbf{0}$ for all $\mathbf{P}$ in $\mathcal{P}$ and all $\mathbf{Q}_1 \neq \mathbf{Q}_2$ in $\mathcal{Q}$.

The structure $\mathcal{Q}$ is *first-order balanced* in relation to $\mathcal{P}$ if (i) holds, and the structure $\mathcal{Q}$ is *orthogonal* in relation to $\mathcal{P}$ if (i) and (ii) hold with each $\lambda_{\mathbf{PQ}}$ equal to either 1 or 0.

Of course, these definitions implicitly depend on the function $f$. It is useful to summarize the set of efficiency factors for the relationship between the two decompositions in the $\mathcal{P} \times \mathcal{Q}$ *efficiency matrix* $\Lambda_{\mathcal{PQ}}$. This matrix $\Lambda_{\mathcal{PQ}}$ is equivalent to the $\mathcal{P} \times \mathcal{Q}$ table introduced in [24] to summarize the efficiency factors in a generally balanced experiment (see also [7], Chapter 7).

THEOREM 4.1 (James and Wilkinson [26]). *Suppose that $\mathbf{Q}$ has first-order balance in relation to $\mathbf{P}$ with efficiency factor $\lambda_{\mathbf{PQ}}$. If $\lambda_{\mathbf{PQ}} \neq 0$, then $\lambda_{\mathbf{PQ}}$ is the unique nonzero eigenvalue of $\mathbf{PQP}$, and so the matrix of orthogonal projection onto $\mathbf{P}(\operatorname{Im}\mathbf{Q})$ is $\lambda_{\mathbf{PQ}}^{-1}\mathbf{PQP}$; if $\lambda_{\mathbf{PQ}} = 0$, then $\operatorname{Im}\mathbf{Q}$ is orthogonal to $\operatorname{Im}\mathbf{P}$.*

DEFINITION 2. If $\mathbf{P}$ and $\mathbf{Q}$ satisfy (4.2) with $\lambda_{\mathbf{PQ}} \neq 0$, then we define $\mathbf{P} \triangleright \mathbf{Q}$ to be $\lambda_{\mathbf{PQ}}^{-1}\mathbf{PQP}$.

As previously noted, $\mathbf{P} \triangleright \mathbf{Q}$ could be called "$\mathbf{P}$ pertaining to $\mathbf{Q}$." More accurately, the source $\operatorname{Im}(\mathbf{P} \triangleright \mathbf{Q})$ is "the part of $\operatorname{Im}\mathbf{P}$ with which $\operatorname{Im}\mathbf{Q}$ is (partly) confounded." Another way to put it is that it is "the part of $\operatorname{Im}\mathbf{P}$ that contains information about $\operatorname{Im}\mathbf{Q}$."

LEMMA 4.1. *If $\lambda_{\mathbf{PQ}} = 1$ then $\mathbf{P} \triangleright \mathbf{Q} = \mathbf{Q}$ and so $\operatorname{Im}\mathbf{Q} \leq \operatorname{Im}\mathbf{P}$.*

PROOF. If $\lambda_{\mathbf{PQ}} = 1$, then $\mathbf{QPQ} = \mathbf{Q}$, and so $(\mathbf{PQ} - \mathbf{Q})'(\mathbf{PQ} - \mathbf{Q}) = (\mathbf{QP} - \mathbf{Q})(\mathbf{PQ} - \mathbf{Q}) = \mathbf{Q} - \mathbf{QPQ} = \mathbf{0}$. Thus $\mathbf{PQ} - \mathbf{Q} = \mathbf{0}$. Hence $\mathbf{PQ} = \mathbf{Q} = \mathbf{Q}' = \mathbf{QP}$, and so $\mathbf{P} \triangleright \mathbf{Q} = \mathbf{PQP} = \mathbf{PQ} = \mathbf{Q}$. Now if $\mathbf{w} \in \operatorname{Im} \mathbf{Q}$, then $\mathbf{w} = \mathbf{Qw}$, and so $\mathbf{Pw} = \mathbf{PQw} = \mathbf{Qw} = \mathbf{w}$: hence $\mathbf{w} \in \operatorname{Im} \mathbf{P}$. □

If (4.3) and (4.4) both hold for a fixed $\mathbf{P}$ and all $\mathbf{Q}$, then the images of the $\mathbf{P} \triangleright \mathbf{Q}$, for $\mathbf{Q}$ in $\mathcal{Q}$ with $\lambda_{\mathbf{PQ}} \neq 0$, are mutually orthogonal subspaces of $\operatorname{Im} \mathbf{P}$. This leaves the orthogonal complement in $\operatorname{Im} \mathbf{P}$ of their sum: this space is just $(\operatorname{Im} \mathbf{P}) \cap V_\Upsilon^\perp$ and is often called the *residual* subspace in $\operatorname{Im} \mathbf{P}$.

DEFINITION 3. Define $\mathbf{P} \vdash \mathcal{Q}$ to be the projector on the residual subspace in $\operatorname{Im} \mathbf{P}$. If all $\mathbf{Q}$ in $\mathcal{Q}$ satisfy (4.2), it is given by

(4.5) $$\mathbf{P} \vdash \mathcal{Q} = \mathbf{P} - \sideset{}{'}\sum_{\mathbf{Q} \in \mathcal{Q}} \mathbf{P} \triangleright \mathbf{Q},$$

where $\sum'_{\mathbf{Q} \in \mathcal{Q}}$ means summation over $\mathbf{Q}$ in $\mathcal{Q}$ with $\lambda_{\mathbf{PQ}} \neq 0$.

In particular, if $\mathbf{PI}_\mathcal{Q} = \mathbf{0}$ then $\mathbf{P} \vdash \mathcal{Q} = \mathbf{P}$, while if $\mathbf{PI}_\mathcal{Q} = \mathbf{P}$ then $\mathbf{P} \vdash \mathcal{Q} = \mathbf{0}$.

As noted above, $\mathbf{P} \vdash \mathcal{Q}$ could be called "$\mathbf{P}$ orthogonal to $\mathcal{Q}$." More accurately, the source $\operatorname{Im}(\mathbf{P} \vdash \mathcal{Q})$ is "the part of $\operatorname{Im} \mathbf{P}$ that is orthogonal to $V_\Upsilon$." Another way to put it is that it is "the part of $\operatorname{Im} \mathbf{P}$ that contains no information about $\operatorname{Im} \mathbf{Q}$ for any $\mathbf{Q}$ in $\mathcal{Q}$."

Finally, to realize our goal of a single decomposition of $V_\Omega$ using $\mathcal{P}$ and $\mathcal{Q}$, we collect together the $\mathbf{P} \triangleright \mathbf{Q}$, for all $\mathbf{P}$ and $\mathbf{Q}$, and $\mathbf{P} \vdash \mathcal{Q}$, for all $\mathbf{P}$, to make the following definition.

DEFINITION 4. If $\mathcal{P}$, $\mathcal{Q}$ are orthogonal decompositions of $V_\Omega$, $V_\Upsilon$ respectively, and a function $f : \Omega \to \Upsilon$ is given such that $\mathcal{Q}$ is structure balanced in relation to $\mathcal{P}$, then the decomposition $\mathcal{P} \triangleright \mathcal{Q}$ of $V_\Omega$ is defined to be

$$\{\mathbf{P} \triangleright \mathbf{Q} : \mathbf{P} \in \mathcal{P}, \mathbf{Q} \in \mathcal{Q}, \lambda_{\mathbf{PQ}} \neq 0\} \cup \{\mathbf{P} \vdash \mathcal{Q} : \mathbf{P} \in \mathcal{P}\}.$$

In [29] it is argued that the analysis of data $\mathbf{y}$ from the experiment begins by projecting $\mathbf{y}$ onto the images of the elements of this decomposition.

LEMMA 4.2. *If $|\Upsilon| = |\Omega|$ and $\mathcal{Q}$ is structure balanced in relation to $\mathcal{P}$, then $\mathcal{P} \triangleright \mathcal{Q} = \mathcal{Q}$.*

PROOF. If $|\Upsilon| = |\Omega|$ then $\mathbf{I}_\mathcal{Q} = \mathbf{I}_\mathcal{P}$. If $\mathcal{Q}$ is structure balanced in relation to $\mathcal{P}$, then equation (4.2) shows that $\lambda_{\mathbf{PQ}} \mathbf{QP} = \mathbf{QPI}_\mathcal{Q} \mathbf{P} = \mathbf{QPI}_\mathcal{P} \mathbf{P} = \mathbf{QP}$, so $\lambda_{\mathbf{PQ}}$ is equal to 0 or 1 for all $\mathbf{P}$ in $\mathcal{P}$ and all $\mathbf{Q}$ in $\mathcal{Q}$. Moreover, $\mathbf{PI}_\mathcal{Q} = \mathbf{PI}_\mathcal{P} = \mathbf{P}$, so there are no residual sources. Hence $\mathcal{Q}$ is a refinement of $\mathcal{P}$ in the sense that each subspace in $\mathcal{P}$ is a direct sum of one or more subspaces in $\mathcal{Q}$. □

In [14], Section 8.4, we commented that the direction of randomization is arbitrary when the two sets of objects have the same size. Lemma 4.2 suggests that, in this case, the unrandomized tier should be the one with less structure, that is, the one placing fewer restrictions on the choice of random permutation.

GenStat [32] accepts from the user the specification of $\mathcal{P}$ and $\mathcal{Q}$ via the BLOCKSTRUCTURE and TREATMENTSTRUCTURE directives, respectively, although each of $\mathcal{P}$ and $\mathcal{Q}$ has only to consist of a set of symmetric idempotents that are first-order balanced in the sense that any two elements from different sets meet condition (4.3). It also allows pseudofactors in its TREATMENTSTRUCTURE directive. Its ANOVA directive uses the algorithm described by Wilkinson [39] and Payne and Wilkinson [34], first to check for first-order balance, and then to decompose $\mathbf{y}$ into its projections by the elements of $\mathcal{P} \rhd \mathcal{Q}$. A similar system, implemented in P-Stat, is described in [22]—the Treatment statement specifies $\mathcal{Q}$, and the Block and Error statements together specify $\mathcal{P}$. In the programs R [36] and S-Plus [25], the decomposition $\mathcal{P} \rhd \mathcal{Q}$ is obtained by specifying $\mathcal{P}$ as the part of the model inside an Error function and $\mathcal{Q}$ as the part outside this function. In these last two packages, there can even be nonorthogonality between the idempotents.

In our examples, we often replace $\mathbf{P}$ and $\mathbf{Q}$ in $\lambda_{\mathbf{PQ}}$ by the names of their corresponding sources. Similarly, in $\mathbf{P} \rhd \mathbf{Q}$ and $\mathbf{P} \vdash \mathcal{Q}$ we replace $\mathbf{P}$ and $\mathbf{Q}$ by the names of their sources to indicate the subspaces $\mathbf{P}(\operatorname{Im}\mathbf{Q})$ and $(\operatorname{Im}\mathbf{P}) \cap V_\Upsilon^\perp$, respectively.

EXAMPLE 2 (Continued). The Hasse diagram derived from the middle panel in Figure 2 is in Figure 6, where Q denotes Squares. There are seven elements of $\mathcal{P}$, one for each source in Figure 6. The sources for the treatments tier are in Figure 4, and so $\mathcal{Q}$ has four elements.

The relationship between $\mathcal{P}$ and $\mathcal{Q}$ is exhibited in the efficiency factors between their elements. The nonzero efficiency factors are $\lambda_{\mathrm{Mean,Mean}} = 1$, $\lambda_{\mathrm{C[Q],T}} = 1/9$, $\lambda_{\mathrm{R\#C[Q],T}} = 8/9$, $\lambda_{\mathrm{H[Q\wedge R\wedge C],M}} = 1$ and $\lambda_{\mathrm{H[Q\wedge R\wedge C],T\#M}} = 1$. These lead to the decomposition in Table 1, which we have already discussed. Efficiency factors are not usually shown in the decomposition table if $\mathcal{Q}$ is orthogonal in relation to $\mathcal{P}$; otherwise, only the nonzero efficiency factors are shown. The decomposition of the vector space indexed by the 48 half-plots is specified by the 11 elements in $\mathcal{P} \rhd \mathcal{Q}$, one for each line in the table:

$$\mathcal{P} \rhd \mathcal{Q} = \left\{ \begin{array}{c} \mathbf{P}_{\mathrm{Mean}} \rhd \mathbf{Q}_{\mathrm{Mean}}, \mathbf{P}_{\mathrm{R}} \vdash \mathcal{Q}, \mathbf{P}_{\mathrm{Q}} \vdash \mathcal{Q}, \mathbf{P}_{\mathrm{Q\#R}} \vdash \mathcal{Q}, \\ \mathbf{P}_{\mathrm{C[Q]}} \rhd \mathbf{Q}_{\mathrm{T}}, \mathbf{P}_{\mathrm{C[Q]}} \vdash \mathcal{Q}, \mathbf{P}_{\mathrm{R\#C[Q]}} \rhd \mathbf{Q}_{\mathrm{T}}, \mathbf{P}_{\mathrm{R\#C[Q]}} \vdash \mathcal{Q}, \\ \mathbf{P}_{\mathrm{H[Q\wedge R\wedge C]}} \rhd \mathbf{Q}_{\mathrm{M}}, \mathbf{P}_{\mathrm{H[Q\wedge R\wedge C]}} \rhd \mathbf{Q}_{\mathrm{T\#M}}, \mathbf{P}_{\mathrm{H[Q\wedge R\wedge C]}} \vdash \mathcal{Q} \end{array} \right\}$$

with

$$\mathbf{P}_{\text{Mean}} \triangleright \mathbf{Q}_{\text{Mean}} = \mathbf{P}_{\text{Mean}} = \mathbf{Q}_{\text{Mean}},$$
$$\mathbf{P}_{\text{R}} \vdash \mathcal{Q} = \mathbf{P}_{\text{R}}, \qquad \mathbf{P}_{\text{Q}} \vdash \mathcal{Q} = \mathbf{P}_{\text{Q}}, \qquad \mathbf{P}_{\text{Q\#R}} \vdash \mathcal{Q} = \mathbf{P}_{\text{Q\#R}},$$
$$\mathbf{P}_{\text{C[Q]}} \triangleright \mathbf{Q}_{\text{T}} = (1/9)^{-1} \mathbf{P}_{\text{C[Q]}} \mathbf{Q}_{\text{T}} \mathbf{P}_{\text{C[Q]}},$$
$$\mathbf{P}_{\text{C[Q]}} \vdash \mathcal{Q} = \mathbf{P}_{\text{C[Q]}} - \mathbf{P}_{\text{C[Q]}} \triangleright \mathbf{Q}_{\text{T}},$$
$$\mathbf{P}_{\text{R\#C[Q]}} \triangleright \mathbf{Q}_{\text{T}} = (8/9)^{-1} \mathbf{P}_{\text{R\#C[Q]}} \mathbf{Q}_{\text{T}} \mathbf{P}_{\text{R\#C[Q]}},$$
$$\mathbf{P}_{\text{R\#C[Q]}} \vdash \mathcal{Q} = \mathbf{P}_{\text{R\#C[Q]}} - \mathbf{P}_{\text{R\#C[Q]}} \triangleright \mathbf{Q}_{\text{T}},$$
$$\mathbf{P}_{\text{H[Q}\wedge\text{R}\wedge\text{C]}} \triangleright \mathbf{Q}_{\text{M}} = \mathbf{Q}_{\text{M}}, \qquad \mathbf{P}_{\text{H[Q}\wedge\text{R}\wedge\text{C]}} \triangleright \mathbf{Q}_{\text{T\#M}} = \mathbf{Q}_{\text{T\#M}} \quad \text{and}$$
$$\mathbf{P}_{\text{H[Q}\wedge\text{R}\wedge\text{C]}} \vdash \mathcal{Q} = \mathbf{P}_{\text{H[Q}\wedge\text{R}\wedge\text{C]}} - \mathbf{Q}_{\text{M}} - \mathbf{Q}_{\text{T\#M}}.$$

**5. Structure for composed randomizations.** For a pair of composed randomizations, there are three sets of objects $\Omega$, $\Upsilon$ and $\Gamma$ with structures $\mathcal{P}$, $\mathcal{Q}$ and $\mathcal{R}$, respectively. While the results apply to quite general structures, those for our examples will be derived from a tier of factors that together uniquely index the set of objects. There is one randomization from $\Gamma$ onto $\Upsilon$, and one from $\Upsilon$ to $\Omega$. Two functions are required to encapsulate the results of these randomizations, say $f:\Omega \to \Upsilon$ and $g:\Upsilon \to \Gamma$. This is represented diagrammatically in Figure 7. Of course, there is a third function implied here, namely the composite function $g \circ f:\Omega \to \Gamma$.

Here we extend the results obtained in Section 4 for two structures to obtain a single decomposition of $V_\Omega$ that combines the three structures in composed randomizations. We will write everything as if it were defined on $\Omega$ rather than on $\Upsilon$ or $\Gamma$. In Section 4, we have already justified writing $V_\Upsilon^f$ as $V_\Upsilon$ and $\mathbf{Q}^f$ as $\mathbf{Q}$ when $\Upsilon$ is randomized to $\Omega$, using the function $f$. When $\Gamma$ is randomized to $\Upsilon$, using the function $g$, and $\Upsilon$ is randomized to $\Omega$, using the function $f$, then the embedding of $V_\Gamma$ in $V_\Omega$ depends on the composite function $g \circ f$ but on no further information about how $V_\Gamma$ is embedded in $V_\Upsilon$. Thus we may identify $V_\Gamma$ with $(V_\Gamma^g)^f$. After these identifications we have $V_\Gamma \leq V_\Upsilon \leq V_\Omega$. Hence $\mathbf{I}_\mathcal{Q}$ is the identity not only for all $\mathbf{Q}$ in $\mathcal{Q}$ but also for all $\mathbf{R}$ in $\mathcal{R}$; that is, $\mathbf{RI}_\mathcal{Q} = \mathbf{I}_\mathcal{Q}\mathbf{R} = \mathbf{R}$ for all $\mathbf{R}$ in $\mathcal{R}$.

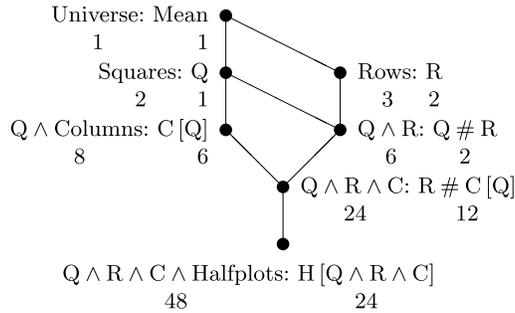

Fig. 6. *Hasse diagram for half-plots in Example 2.*

THEOREM 5.1. *Suppose that corresponding to the sets of objects $\Omega$, $\Upsilon$, $\Gamma$ there are sets $\mathcal{P}$, $\mathcal{Q}$, $\mathcal{R}$ of $\Omega \times \Omega$ matrices, specifying orthogonal decompositions of the subspaces $V_\Omega$, $V_\Upsilon$, $V_\Gamma$, respectively, with $V_\Gamma \leq V_\Upsilon \leq V_\Omega$; that $\mathcal{Q}$ is structure balanced in relation to $\mathcal{P}$ with efficiency matrix $\Lambda_{\mathcal{PQ}}$; and that $\mathcal{R}$ is structure balanced in relation to $\mathcal{Q}$ with efficiency matrix $\Lambda_{\mathcal{QR}}$. Then:*

(a) *$\mathcal{R}$ is structure balanced in relation to $\mathcal{P}$ with efficiency matrix $\Lambda_{\mathcal{PR}} = \Lambda_{\mathcal{PQ}}\Lambda_{\mathcal{QR}}$;*
(b) *$\mathcal{R}$ is structure balanced in relation to $\mathcal{P} \triangleright \mathcal{Q}$ with efficiency matrix $\Lambda_{\mathcal{P} \triangleright \mathcal{Q}, \mathcal{R}}$ whose elements are the efficiency factors $\lambda_{\mathbf{P} \triangleright \mathbf{Q}, \mathbf{R}} = \lambda_{\mathbf{PQ}}\lambda_{\mathbf{QR}}$ and $\lambda_{\mathbf{P} \vdash \mathcal{Q}, \mathbf{R}} = 0$;*
(c) *$\mathcal{Q} \triangleright \mathcal{R}$ is structure balanced in relation to $\mathcal{P}$ with efficiency matrix $\Lambda_{\mathcal{P}, \mathcal{Q} \triangleright \mathcal{R}}$ whose elements are the efficiency factors $\lambda_{\mathbf{P}, \mathbf{Q} \triangleright \mathbf{R}} = \lambda_{\mathbf{PQ}}$ and $\lambda_{\mathbf{P}, \mathbf{Q} \vdash \mathcal{R}} = \lambda_{\mathbf{PQ}}$;*
(d) *$(\mathcal{P} \triangleright \mathcal{Q}) \triangleright \mathcal{R} = \mathcal{P} \triangleright (\mathcal{Q} \triangleright \mathcal{R}) =*
$$\{(\mathbf{P} \triangleright \mathbf{Q}) \triangleright \mathbf{R} : \mathbf{P} \in \mathcal{P}, \mathbf{Q} \in \mathcal{Q}, \mathbf{R} \in \mathcal{R}, \lambda_{\mathbf{PQ}} \neq 0, \lambda_{\mathbf{QR}} \neq 0\}$$
$$\cup \{(\mathbf{P} \triangleright \mathbf{Q}) \vdash \mathcal{R} : \mathbf{P} \in \mathcal{P}, \mathbf{Q} \in \mathcal{Q}, \lambda_{\mathbf{PQ}} \neq 0\} \cup \{\mathbf{P} \vdash \mathcal{Q} : \mathbf{P} \in \mathcal{P}\}.$$

PROOF. (a) We need to show that $\mathbf{RPR} = \sum_{\mathbf{Q} \in \mathcal{Q}} \lambda_{\mathbf{PQ}}\lambda_{\mathbf{QR}}\mathbf{R}$ for $\mathbf{P}$ in $\mathcal{P}$ and $\mathbf{R}$ in $\mathcal{R}$, and that $\mathbf{R}_1\mathbf{PR}_2 = \mathbf{0}$ for $\mathbf{R}_1$, $\mathbf{R}_2$ in $\mathcal{R}$ with $\mathbf{R}_1 \neq \mathbf{R}_2$.

Since $\sum_{\mathbf{Q} \in \mathcal{Q}} \mathbf{Q} = \mathbf{I}_\mathcal{Q}$, we have $\mathbf{I}_\mathcal{Q}\mathbf{PI}_\mathcal{Q} = \sum_{\mathbf{Q} \in \mathcal{Q}} \mathbf{QPI}_\mathcal{Q} = \sum_{\mathbf{Q} \in \mathcal{Q}} \lambda_{\mathbf{PQ}}\mathbf{Q}$. Also $\mathbf{R} = \mathbf{RI}_\mathcal{Q}$. Hence
$$\mathbf{RPR} = \mathbf{RI}_\mathcal{Q}\mathbf{PI}_\mathcal{Q}\mathbf{R} = \sum_\mathbf{Q} \mathbf{R}\lambda_{\mathbf{PQ}}\mathbf{QR} = \sum_\mathbf{Q} \lambda_{\mathbf{PQ}}\lambda_{\mathbf{QR}}\mathbf{R}$$
while $\mathbf{R}_1\mathbf{PR}_2 = \sum_\mathbf{Q} \lambda_{\mathbf{PQ}}\mathbf{R}_1\mathbf{QR}_2 = \mathbf{0}$.

(b) For $\mathcal{R}$ to be structure balanced in relation to $\mathcal{P} \triangleright \mathcal{Q}$ with the given efficiency factors, we require the following conditions for all $\mathbf{P}$ in $\mathcal{P}$, $\mathbf{Q}$ in $\mathcal{Q}$, and $\mathbf{R}$, $\mathbf{R}_i$, $\mathbf{R}_j$ in $\mathcal{R}$: (i) $\mathbf{R}(\mathbf{P} \triangleright \mathbf{Q})\mathbf{R} = \lambda_{\mathbf{PQ}}\lambda_{\mathbf{QR}}\mathbf{R}$; (ii) $\mathbf{R}_i(\mathbf{P} \triangleright \mathbf{Q})\mathbf{R}_j = \mathbf{0}$ if $\mathbf{R}_i \neq \mathbf{R}_j$; (iii) $\mathbf{R}_i(\mathbf{P} \vdash \mathcal{Q})\mathbf{R}_j = \mathbf{0}$.

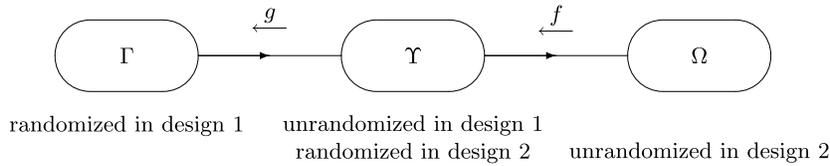

FIG. 7. *Diagram of an experiment with two composed randomizations.*

Equation (4.2) shows that $\mathbf{I}_\mathcal{Q}(\mathbf{P} \triangleright \mathbf{Q})\mathbf{I}_\mathcal{Q} = \mathbf{I}_\mathcal{Q}\lambda_{\mathbf{PQ}}^{-1}\mathbf{PQPI}_\mathcal{Q} = \lambda_{\mathbf{PQ}}\mathbf{Q}$. Therefore $\mathbf{R}(\mathbf{P} \triangleright \mathbf{Q})\mathbf{R} = \mathbf{RI}_\mathcal{Q}(\mathbf{P} \triangleright \mathbf{Q})\mathbf{I}_\mathcal{Q}\mathbf{R} = \lambda_{\mathbf{PQ}}\mathbf{RQR} = \lambda_{\mathbf{PQ}}\lambda_{\mathbf{QR}}\mathbf{R}$ while $\mathbf{R}_1(\mathbf{P} \triangleright \mathbf{Q})\mathbf{R}_2 = \lambda_{\mathbf{PQ}}\mathbf{R}_1\mathbf{QR}_2 = \mathbf{0}$ if $\mathbf{R}_1 \neq \mathbf{R}_2$. This proves (i) and (ii).

From equation (4.5), $\mathbf{R}_i\mathbf{PR}_j - \mathbf{R}_i(\mathbf{P} \vdash \mathcal{Q})\mathbf{R}_j = \sum'_\mathbf{Q} \mathbf{R}_i(\mathbf{P} \triangleright \mathbf{Q})\mathbf{R}_j$. Hence part (a) shows that $\mathbf{R}_i(\mathbf{P} \vdash \mathcal{Q})\mathbf{R}_j = \mathbf{0}$ if $\mathbf{R}_i \neq \mathbf{R}_j$ while $\mathbf{R}(\mathbf{P} \vdash \mathcal{Q})\mathbf{R} = \mathbf{RPR} - \sum_\mathbf{Q} \lambda_{\mathbf{PQ}}\lambda_{\mathbf{QR}}\mathbf{R} = \mathbf{0}$. This proves (iii).

(c) If $\mathbf{M}$ is in $\mathcal{Q} \triangleright \mathcal{R}$ then $\mathbf{I}_\mathcal{Q}\mathbf{M} = \mathbf{M} = \mathbf{MI}_\mathcal{Q}$, so $\mathbf{M}_i\mathbf{PM}_j = \mathbf{M}_i\mathbf{I}_\mathcal{Q}\mathbf{PI}_\mathcal{Q}\mathbf{M}_j = \mathbf{M}_i\sum_\mathbf{Q}\lambda_{\mathbf{PQ}}\mathbf{QM}_j$, which is zero unless $\mathbf{M}_i$ and $\mathbf{M}_j$ have the same $\mathbf{Q}$. If $\mathbf{M}_i$ and $\mathbf{M}_j$ have the same $\mathbf{Q}$ then $\mathbf{M}_i\mathbf{PM}_j = \lambda_{\mathbf{PQ}}\mathbf{M}_i\mathbf{QM}_j = \lambda_{\mathbf{PQ}}\mathbf{M}_i\mathbf{M}_j = \delta_{ij}\lambda_{\mathbf{PQ}}\mathbf{M}_i\mathbf{M}_j$. Hence $\mathcal{Q} \triangleright \mathcal{R}$ is structure balanced in relation to $\mathcal{P}$ with $\lambda_{\mathbf{P},\mathbf{Q}\triangleright\mathbf{R}} = \lambda_{\mathbf{PQ}}$ and $\lambda_{\mathbf{P},\mathbf{Q}\vdash\mathcal{R}} = \lambda_{\mathbf{PQ}}$.

(d) By (b), we may apply Theorem 4.1 and Definition 2 to the pair of structures $\mathcal{P} \triangleright \mathcal{Q}$ and $\mathcal{R}$ to obtain the matrix $(\mathbf{P} \triangleright \mathbf{Q}) \triangleright \mathbf{R}$ of orthogonal projection onto $(\mathbf{P} \triangleright \mathbf{Q})(\operatorname{Im} \mathbf{R})$:

$$(\mathbf{P} \triangleright \mathbf{Q}) \triangleright \mathbf{R} = \lambda_{\mathbf{PQ}}^{-1}\lambda_{\mathbf{QR}}^{-1}(\mathbf{P} \triangleright \mathbf{Q})\mathbf{R}(\mathbf{P} \triangleright \mathbf{Q}) = \lambda_{\mathbf{PQ}}^{-1}\lambda_{\mathbf{QR}}^{-1}(\lambda_{\mathbf{PQ}}^{-1}\mathbf{PQP})\mathbf{R}(\lambda_{\mathbf{PQ}}^{-1}\mathbf{PQP}).$$

Moreover, equation (4.5) gives

$$(\mathbf{P} \triangleright \mathbf{Q}) \vdash \mathcal{R} = \mathbf{P} \triangleright \mathbf{Q} - \sum_{\substack{\mathbf{R} \in \mathcal{R} \\ \lambda_{\mathbf{QR}} \neq 0}} (\mathbf{P} \triangleright \mathbf{Q}) \triangleright \mathbf{R}.$$

Since each $\mathbf{R}$ is orthogonal to $\mathbf{P} \vdash \mathcal{Q}$, we have $(\mathbf{P} \vdash \mathcal{Q}) \vdash \mathcal{R} = \mathbf{P} \vdash \mathcal{Q}$. Thus Definition 4 gives the decomposition $(\mathcal{P} \triangleright \mathcal{Q}) \triangleright \mathcal{R}$.

Similarly, for the pair of structures $\mathcal{P}$ and $\mathcal{Q} \triangleright \mathcal{R}$, part (c) gives

$$\mathbf{P} \triangleright (\mathbf{Q} \triangleright \mathbf{R}) = \lambda_{\mathbf{PQ}}^{-1}\mathbf{P}(\mathbf{Q} \triangleright \mathbf{R})\mathbf{P} = \lambda_{\mathbf{PQ}}^{-1}\mathbf{P}(\lambda_{\mathbf{QR}}^{-1}\mathbf{QRQ})\mathbf{P}$$

and

$$\mathbf{P} \triangleright (\mathbf{Q} \vdash \mathcal{R}) = \lambda_{\mathbf{PQ}}^{-1}\mathbf{P}\Big(\mathbf{Q} - \sum_{\substack{\mathbf{R} \in \mathcal{R} \\ \lambda_{\mathbf{QR}} \neq 0}} \mathbf{Q} \triangleright \mathbf{R}\Big)\mathbf{P}.$$

Then, by applying equation (4.5) twice,

$$\mathbf{P} \vdash (\mathcal{Q} \triangleright \mathcal{R}) = \mathbf{P} - \sum_{\mathbf{R} \in \mathcal{R}}\sum_{\mathbf{Q} \in \mathcal{Q}} \mathbf{P} \triangleright (\mathbf{Q} \triangleright \mathbf{R}) - \sum_{\mathbf{Q} \in \mathcal{Q}} \mathbf{P} \triangleright (\mathbf{Q} \vdash \mathcal{R})$$

$$= \mathbf{P} - \sum_{\mathbf{Q} \in \mathcal{Q}} \mathbf{P} \triangleright \mathbf{Q} = \mathbf{P} \vdash \mathcal{Q}.$$

Hence Definition 4 gives $\mathcal{P} \triangleright (\mathcal{Q} \triangleright \mathcal{R}) =$

(5.1) $$\{\mathbf{P} \triangleright (\mathbf{Q} \triangleright \mathbf{R}) : \mathbf{P} \in \mathcal{P}, \mathbf{Q} \in \mathcal{Q}, \mathbf{R} \in \mathcal{R}, \lambda_{\mathbf{PQ}} \neq 0, \lambda_{\mathbf{QR}} \neq 0\}$$
$$\cup \{\mathbf{P} \triangleright (\mathbf{Q} \vdash \mathcal{R}) : \mathbf{P} \in \mathcal{P}, \mathbf{Q} \in \mathcal{Q}, \lambda_{\mathbf{PQ}} \neq 0\} \cup \{\mathbf{P} \vdash \mathcal{Q} : \mathbf{P} \in \mathcal{P}\}.$$

It now suffices to show that $(\mathbf{P} \triangleright \mathbf{Q}) \triangleright \mathbf{R} = \mathbf{P} \triangleright (\mathbf{Q} \triangleright \mathbf{R})$ and $(\mathbf{P} \triangleright \mathbf{Q}) \vdash \mathcal{R} = \mathbf{P} \triangleright (\mathbf{Q} \vdash \mathcal{R})$ for all $\mathbf{P}$ in $\mathcal{P}$, $\mathbf{Q}$ in $\mathcal{Q}$ and $\mathbf{R}$ in $\mathcal{R}$. Now, $\mathbf{QPR} = \mathbf{QPI}_\mathcal{Q}\mathbf{R} = \lambda_{\mathbf{PQ}}\mathbf{QR}$, so

$$(\mathbf{P} \triangleright \mathbf{Q}) \triangleright \mathbf{R} = \lambda_{\mathbf{PQ}}^{-3}\lambda_{\mathbf{QR}}^{-1}\mathbf{PQPRPQP} = \lambda_{\mathbf{PQ}}^{-1}\lambda_{\mathbf{QR}}^{-1}\mathbf{PQRQP} = \mathbf{P} \triangleright (\mathbf{Q} \triangleright \mathbf{R}).$$

Also,

$$(\mathbf{P} \triangleright \mathbf{Q}) \vdash \mathcal{R} = \mathbf{P} \triangleright \mathbf{Q} - \sum_{\mathbf{R}\in\mathcal{R}}(\mathbf{P} \triangleright \mathbf{Q}) \triangleright \mathbf{R} = \mathbf{P} \triangleright \mathbf{Q} - \sum_{\mathbf{R}\in\mathcal{R}}\mathbf{P} \triangleright (\mathbf{Q} \triangleright \mathbf{R})$$

$$= \lambda_{\mathbf{PQ}}^{-1}\mathbf{P}\Big(\mathbf{Q} - {\sum_{\mathbf{R}\in\mathcal{R}}}'\mathbf{Q} \triangleright \mathbf{R}\Big)\mathbf{P} = \mathbf{P} \triangleright (\mathbf{Q} \vdash \mathcal{R}). \qquad \square$$

As noted in Section 4, if $\mathbf{PI}_\mathcal{Q} = \mathbf{0}$ then $\mathbf{P} \vdash \mathcal{Q} = \mathbf{P}$ while if $\mathbf{PI}_\mathcal{Q} = \mathbf{P}$ then $\mathbf{P} \vdash \mathcal{Q} = \mathbf{0}$. Analogously, if $(\mathbf{P} \triangleright \mathbf{Q})\mathbf{I}_\mathcal{R} = \mathbf{0}$ then $(\mathbf{P} \triangleright \mathbf{Q}) \vdash \mathcal{R} = \mathbf{P} \triangleright \mathbf{Q}$ while if $(\mathbf{P} \triangleright \mathbf{Q})\mathbf{I}_\mathcal{R} = \mathbf{P} \triangleright \mathbf{Q}$ then $(\mathbf{P} \triangleright \mathbf{Q}) \vdash \mathcal{R} = \mathbf{0}$.

The algorithm in [39] for analysing two-tiered experiments obtains the decomposition of the data vector $\mathbf{y}$ into orthogonal vectors of the form $\mathbf{Py}$ according to the decomposition $\mathcal{P}$ and then further decomposes the result according to $\mathcal{Q}$, to yield the combined decomposition $\mathcal{P} \triangleright \mathcal{Q}$ of $\mathbf{y}$. If $\mathcal{R}$ is structure balanced in relation to $\mathcal{P} \triangleright \mathcal{Q}$, then we can run the algorithm again to obtain the decomposition $(\mathcal{P} \triangleright \mathcal{Q}) \triangleright \mathcal{R}$ of $\mathbf{y}$. However, consideration of the order in which the randomizations have been performed suggests that $\mathcal{P}$ should be decomposed according to the joint decomposition $\mathcal{Q} \triangleright \mathcal{R}$, to yield the decomposition $\mathcal{P} \triangleright (\mathcal{Q} \triangleright \mathcal{R})$ of $\mathbf{y}$. Theorem 5.1(d) shows that the decomposition may be performed in either order.

EXAMPLE 1 (Continued). The Hasse diagrams for the treatments and meatloaves tiers are like those in Figures 4 and 3, respectively, while Figure 8 gives the Hasse diagram for the tastings tier. Let $\mathcal{P}$, $\mathcal{Q}$ and $\mathcal{R}$ be the structures corresponding to tastings, meatloaves and treatments, respectively. Both designs are orthogonal. Table 2 gives the decomposition table corresponding to the joint decomposition $(\mathcal{P} \triangleright \mathcal{Q}) \triangleright \mathcal{R}$, or, equivalently, $\mathcal{P} \triangleright (\mathcal{Q} \triangleright \mathcal{R})$. For the former decomposition, first obtain $\mathcal{P} \triangleright \mathcal{Q}$, shown in the two left columns. Then $\mathcal{P} \triangleright \mathcal{Q}$ is decomposed by $\mathcal{R}$.

In this case, every nonzero entry in $\Lambda_{\mathcal{PQ}}$ and $\Lambda_{\mathcal{QR}}$ is equal to one. There is exactly one nonzero entry in each column of $\Lambda_{\mathcal{QR}}$, so the nonzero entries in the product $\Lambda_{\mathcal{PQ}}\Lambda_{\mathcal{QR}}$ are also all equal to one. Theorem 5.1(a) states that $\Lambda_{\mathcal{PR}} = \Lambda_{\mathcal{PQ}}\Lambda_{\mathcal{QR}}$, so we see that, if the two component designs are orthogonal, then so is the composite design, which ignores the middle tier. In other words, $\mathcal{R}$ is orthogonal in relation to $\mathcal{P}$.

Moreover, the nonzero products of entries of $\Lambda_{\mathcal{PQ}}$ with entries of $\Lambda_{QR}$ are also all equal to one. Hence Theorem 5.1(b) shows that $\Lambda_{\mathcal{P}\triangleright\mathcal{Q},\mathcal{R}}$ also has all nonzero elements equal to one. In other words, $\mathcal{R}$ is orthogonal in relation to $\mathcal{P} \triangleright \mathcal{Q}$, and thus there is no need to show efficiency factors in the decomposition table.

Table 2
*Decomposition table for Example 1*

| tastings tier | | meatloaves tier | | treatments tier | |
|---|---|---|---|---|---|
| **source** | **d.f.** | **source** | **d.f.** | **source** | **d.f.** |
| Mean | 1 | Mean | 1 | Mean | 1 |
| Sessions | 2 | Blocks | 2 | | |
| Panellists[S] | 33 | | | | |
| Time-orders[S] | 15 | | | | |
| P#T[S] | 165 | Meatloaves[B] | 15 | Rosemary | 1 |
| | | | | Irradiation | 2 |
| | | | | Rosemary#Irradiation | 2 |
| | | | | Residual | 10 |
| | | Residual | 150 | | |

The full decomposition contains nine elements, one for each line in the table:

$$(\mathcal{P} \triangleright \mathcal{Q}) \triangleright \mathcal{R} = \left\{ \begin{array}{c} (\mathbf{P}_{\text{Mean}} \triangleright \mathbf{Q}_{\text{Mean}}) \triangleright \mathbf{R}_{\text{Mean}}, \mathbf{P}_S \triangleright \mathbf{Q}_B, \mathbf{P}_{P[S]}, \mathbf{P}_{T[S]}, \\ (\mathbf{P}_{P\#T[S]} \triangleright \mathbf{Q}_{M[B]}) \triangleright \mathbf{R}_R, (\mathbf{P}_{P\#T[S]} \triangleright \mathbf{Q}_{M[B]}) \triangleright \mathbf{R}_I, \\ (\mathbf{P}_{P\#T[S]} \triangleright \mathbf{Q}_{M[B]}) \triangleright \mathbf{R}_{R\#I}, (\mathbf{P}_{P\#T[S]} \triangleright \mathbf{Q}_{M[B]}) \vdash \mathcal{R}, \\ \mathbf{P}_{P\#T[S]} \vdash \mathcal{Q} \end{array} \right\}$$

with

$$(\mathbf{P}_{\text{Mean}} \triangleright \mathbf{Q}_{\text{Mean}}) \triangleright \mathbf{R}_{\text{Mean}} = \mathbf{P}_{\text{Mean}} = \mathbf{Q}_{\text{Mean}} = \mathbf{R}_{\text{Mean}},$$

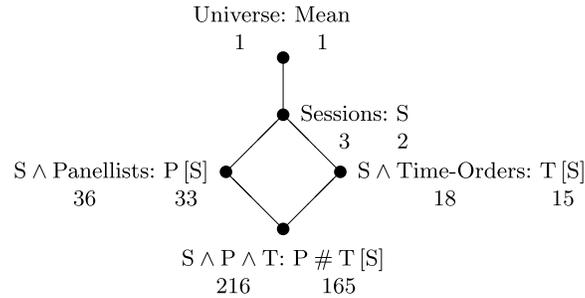

Fig. 8. *Hasse diagram for tastings in Example 1.*

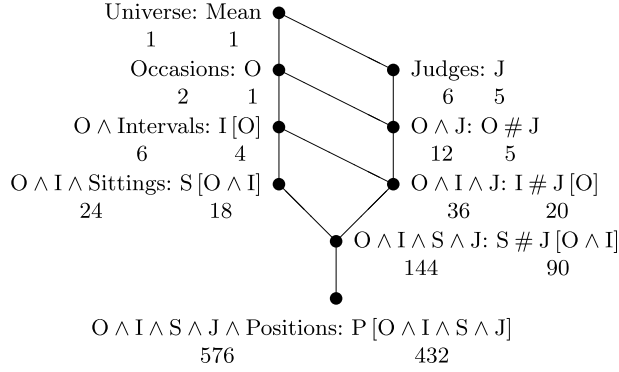

Fig. 9. *Hasse diagram for evaluations in Example 2.*

$$\mathbf{P}_S \triangleright \mathbf{Q}_B = \mathbf{P}_S = \mathbf{Q}_B = (\mathbf{P}_S \triangleright \mathbf{Q}_B) \vdash \mathcal{R},$$
$$\mathbf{P}_{P[S]} = \mathbf{P}_{P[S]} \vdash \mathcal{Q},$$
$$\mathbf{P}_{P[T]} = \mathbf{P}_{P[T]} \vdash \mathcal{Q}, \qquad \mathbf{P}_{P\#T[S]} \triangleright \mathbf{Q}_{M[B]} = \mathbf{Q}_{M[B]},$$
$$(\mathbf{P}_{P\#T[S]} \triangleright \mathbf{Q}_{M[B]}) \triangleright \mathbf{R}_R = \mathbf{R}_R, \qquad (\mathbf{P}_{P\#T[S]} \triangleright \mathbf{Q}_{M[B]}) \triangleright \mathbf{R}_I = \mathbf{R}_I,$$
$$(\mathbf{P}_{P\#T[S]} \triangleright \mathbf{Q}_{M[B]}) \triangleright \mathbf{R}_{R\#I} = \mathbf{R}_{R\#I},$$
$$(\mathbf{P}_{P\#T[S]} \triangleright \mathbf{Q}_{M[B]}) \vdash \mathcal{R} = \mathbf{Q}_{M[B]} \vdash \mathcal{R} = \mathbf{Q}_{M[B]} - \mathbf{R}_R - \mathbf{R}_I - \mathbf{R}_{R\#I} \quad \text{and}$$
$$\mathbf{P}_{P\#T[S]} \vdash \mathcal{Q} = \mathbf{P}_{P\#T[S]} - \mathbf{Q}_{M[B]}.$$

Note the occurrence of the three different types of elements in the decomposition: $(\mathbf{P} \triangleright \mathbf{Q}) \triangleright \mathbf{R}$, $(\mathbf{P} \triangleright \mathbf{Q}) \vdash \mathcal{R}$ and $\mathbf{P} \vdash \mathcal{Q}$. As demonstrated here, for an orthogonal experiment some projection operators for the decomposition are just elements of $\mathcal{P}$, $\mathcal{Q}$ and $\mathcal{R}$ while each of the others is obtained by subtracting some of these elements from another one.

EXAMPLE 2 (Continued). The panel on the right of Figure 2 gives the Hasse diagram in Figure 9. This gives the sources in the left column of Table 3.

Let $\mathcal{P}$, $\mathcal{Q}$ and $\mathcal{R}$ be the structures corresponding to evaluations, half-plots and treatments, respectively. As $\mathcal{Q}$ is structure balanced in relation to $\mathcal{P}$, and $\mathcal{R}$ is structure balanced in relation to $\mathcal{Q}$, Theorem 5.1(a) implies that $\mathcal{R}$ is also structure balanced in relation to $\mathcal{P}$ and that the matrix $\Lambda_{\mathcal{P}\mathcal{R}}$ of efficiencies for treatments in relation to evaluations is $\Lambda_{\mathcal{P}\mathcal{Q}}\Lambda_{\mathcal{Q}\mathcal{R}}$. Computing this product shows that Trellis is estimated from $S[O \wedge I]$ with efficiency $1/27$ and from $S\#J[O \wedge I]$ with efficiency $26/27$.

Table 3 gives the joint decomposition $(\mathcal{P} \triangleright \mathcal{Q}) \triangleright \mathcal{R}$, which contains 20 elements. For this decomposition, first obtain $\mathcal{P} \triangleright \mathcal{Q}$, shown in the left two columns, with the nonzero entries of $\Lambda_{\mathcal{P}\mathcal{Q}}$ included in the middle column.

DECOMPOSITION TABLES I. A CHAIN OF RANDOMIZATIONS        23TABLE 3
*Decomposition table for Example 2 (O = Occasions, I = Intervals, S = Sittings, J = Judges, P = Positions, Q = Squares, C = Columns, R = Rows, H = Halfplots, T = Trellis, M = Methods)*

| evaluations tier | | half-plots tier | | | treatments tier | | |
|---|---|---|---|---|---|---|---|
| **source** | **d.f.** | **eff.** | **source** | **d.f.** | **eff.** | **source** | **d.f.** |
| Mean | 1 | 1 | Mean | 1 | 1 | Mean | 1 |
| O | 1 | 1 | Q | 1 | | | |
| I[O] | 4 | | | | | | |
| S[O ∧ I] | 18 | $\frac{1}{3}$ | C[Q] | 6 | $\frac{1}{27}$ | T | 3 |
| | | | | | | Residual | 3 |
| | | | Residual | 12 | | | |
| J | 5 | | | | | | |
| O#J | 5 | | | | | | |
| I#J[O] | 20 | 1 | R | 2 | | | |
| | | 1 | Q#R | 2 | | | |
| | | | Residual | 16 | | | |
| S#J[O ∧ I] | 90 | $\frac{2}{3}$ | C[Q] | 6 | $\frac{2}{27}$ | T | 3 |
| | | | | | | Residual | 3 |
| | | 1 | R#C[Q] | 12 | $\frac{8}{9}$ | T | 3 |
| | | | | | | Residual | 9 |
| | | | Residual | 72 | | | |
| P[O ∧ I ∧ S ∧ J] | 432 | 1 | H[Q ∧ R ∧ C] | 24 | 1 | M | 1 |
| | | | | | 1 | T#M | 3 |
| | | | | | | Residual | 20 |
| | | | Residual | 408 | | | |

Then $\mathcal{P} \triangleright \mathcal{Q}$ is decomposed by $\mathcal{R}$; the properties of this decomposition are summarized in $\Lambda_{\mathcal{P} \triangleright \mathcal{Q}, \mathcal{R}}$. As Theorem 5.1(b) states, the entries of this matrix are the products of the entries of $\Lambda_{\mathcal{PQ}}$ with those of $\Lambda_{\mathcal{QR}}$. For example,

$$\lambda_{\text{S}[O \wedge I] \triangleright \text{C}[Q], \text{T}} = \lambda_{\text{S}[O \wedge I], \text{C}[Q]} \lambda_{\text{C}[Q], \text{T}} = \tfrac{1}{3} \times \tfrac{1}{9} = \tfrac{1}{27}$$

where Q denotes Squares. These nonzero entries are given in the efficiency column of the treatments tier in Table 3. The two elements of $(\mathcal{P} \triangleright \mathcal{Q}) \triangleright \mathcal{R}$ that correspond to the fourth and sixth lines of sources in Table 3 are:

$$(\mathbf{P}_{\text{S}[O \wedge I]} \triangleright \mathbf{Q}_{\text{C}[Q]}) \triangleright \mathbf{R}_{\text{T}} = (\tfrac{1}{3} \times \tfrac{1}{9})^{-1} \mathbf{P}_{\text{S}[O \wedge I]} \mathbf{Q}_{\text{C}[Q]} \mathbf{R}_{\text{T}} \mathbf{Q}_{\text{C}[Q]} \mathbf{P}_{\text{S}[O \wedge I]}$$

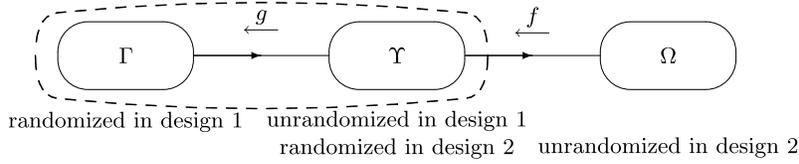

Fig. 10. *Diagram of an experiment with two randomized-inclusive randomizations.*

and

$$\mathbf{P}_{S[O\wedge I]} \vdash \mathcal{Q} = \mathbf{P}_{S[O\wedge I]} - \mathbf{P}_{S[O\wedge I]} \triangleright \mathbf{Q}_{C[Q]}$$
$$= \mathbf{P}_{S[O\wedge I]} - (\tfrac{1}{3})^{-1}\mathbf{P}_{S[O\wedge I]}\mathbf{Q}_{C[Q]}\mathbf{P}_{S[O\wedge I]}.$$

A full analysis-of-variance table is in [15].

**6. Structure for randomized-inclusive randomizations.** The set-up for an experiment with a pair of randomized-inclusive randomizations is the same as for composed randomizations, except that the randomization of $\Upsilon$ to $\Omega$ uses a design whose properties depend on the outcome of the first randomization. The order of the randomizations is prescribed, because the outcome of the first must be known before the second can be done. Two functions are required to encapsulate the results of these randomizations, say $g:\Upsilon \to \Gamma$ and $f:\Omega \to \Upsilon$. This is shown diagrammatically in Figure 10. Note the similarity to Figure 7 for composed randomizations.

There are structures $\mathcal{P}$, $\mathcal{Q}_1$ and $\mathcal{R}_1$ on $\Omega$, $\Upsilon$ and $\Gamma$, respectively. We assume that design 1 is randomized using the group of permutations of $\Upsilon$ corresponding to structure $\mathcal{Q}_1$ and that $\mathcal{R}_1$ is structure balanced in relation to $\mathcal{Q}_1$. For randomized-inclusive randomizations, the second randomization starts with a systematic plan using information from both the previous tiers, and uses the group of permutations of $\Omega$ corresponding to $\mathcal{P}$. However, unlike the situation for composed randomizations, $\mathcal{Q}_1$ is not structure balanced in relation to $\mathcal{P}$.

In simple cases, $\mathcal{Q}_1 \triangleright \mathcal{R}_1$ is structure balanced in relation to $\mathcal{P}$, so the decomposition $\mathcal{P} \triangleright (\mathcal{Q}_1 \triangleright \mathcal{R}_1)$ is appropriate, and the decomposition must proceed from right to left. Example 1 of Wood, Williams and Speed [40] is of this form, as is the randomization in [14], Figure 17. In such cases, the second randomization is done using pseudofactors on $\Upsilon$ of the form $F \circ g$ for some generalized factors $F$ on $\Gamma$: these pseudofactors refine $\mathcal{Q}_1$ to $\mathcal{Q}_1 \triangleright \mathcal{R}_1$. In practice we put $\mathcal{Q} = \mathcal{Q}_1 \triangleright \mathcal{R}_1$ and regard $\mathcal{Q}$ as the appropriate structure on $\Upsilon$. Sometimes the pseudofactors needed on $\Upsilon$ for design 2 are defined by different generalized factors on $\Gamma$ in each replicate. In these cases, $\mathcal{Q}_1$ must be refined to $\mathcal{Q}_2$ by these extra pseudofactors on $\Upsilon$, and then in practice we put $\mathcal{Q} = \mathcal{Q}_2$.

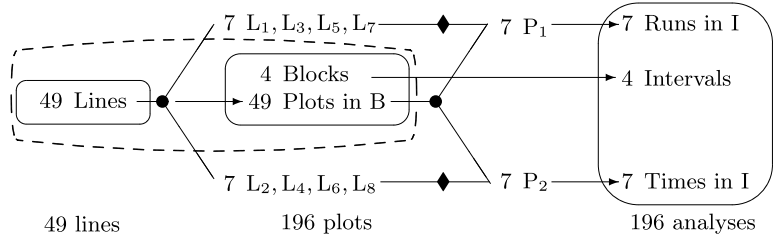

Fig. 11. *Randomized-inclusive randomizations in Example 5: lines are randomized to plots, which are in turn randomized to analyses;* B *denotes Blocks,* I *denotes Intervals;* $L_1, \ldots, L_8$ *are mutually orthogonal pseudofactors for Lines;* $P_1$ *and* $P_2$ *are pseudofactors for Plots, determined from different Lines pseudofactors in different blocks.*

In more complicated cases, design 2 uses pseudofactors on $\Gamma$ that were not used in design 1. These refine $\mathcal{R}_1$ to $\mathcal{R}_2$, and must be used to define pseudofactors on $\Upsilon$ that refine $\mathcal{Q}_1$ to $\mathcal{Q}_2$. Then we put $\mathcal{Q} = \mathcal{Q}_2$, as in Example 5.

Writing $\mathcal{R}$ for $\mathcal{R}_1$ or $\mathcal{R}_2$ as appropriate, in all three cases we have $\mathcal{R}$ structure balanced in relation to $\mathcal{Q}$ and $\mathcal{Q}$ structure balanced in relation to $\mathcal{P}$. Hence Theorem 5.1 shows that $\mathcal{R}$ is structure balanced in relation to $\mathcal{P}$ and that $\Lambda_{\mathcal{PR}} = \Lambda_{\mathcal{PQ}} \Lambda_{\mathcal{QR}}$. The orthogonal decomposition of $V_\Omega$ resulting from the combination of the structures on $\mathcal{P}$, $\mathcal{Q}$ and $\mathcal{R}$ is given by either $(\mathcal{P} \triangleright \mathcal{Q}) \triangleright \mathcal{R}$ or $\mathcal{P} \triangleright (\mathcal{Q} \triangleright \mathcal{R})$, which by Theorem 5.1(d) are equal.

Example 5 (A two-phase wheat variety trial). Example 9 in [14] describes an experiment that consists of a field phase and a laboratory phase. In the field phase, 49 lines of wheat are investigated using a randomized complete-block design with four blocks. The produce of each plot is analysed using a gas chromatograph in which seven samples can be processed per run. In the laboratory phase, blocks are randomized to four intervals. In each interval, there are seven runs, in which samples are processed at seven consecutive times. The design allocates plots to analyses in such a way that the runs and times form a $7 \times 7$ balanced lattice square on the lines. The sets of objects for this experiment are analyses, plots and lines. Pseudofactors are introduced for lines and plots in order to define the design of the second phase. See Figure 11.

Figure 12 shows the Hasse diagrams when pseudofactors are included. The sources for the structure $\mathcal{P}$ are in the top, left-hand diagram. Without pseudofactors, the structure $\mathcal{Q}_1$ on plots has sources Mean, Blocks and Plots[B]. With pseudofactors, the sources corresponding to $\mathcal{Q}_2$ are those in the top, right-hand diagram. Note the source Plots[B]$_\vdash$, which consists of all contrasts for Plots orthogonal not only to Blocks, as for Plots[B], but also to the pseudofactors $P_1$ and $P_2$. Similarly, the natural structure $\mathcal{R}_1$ on lines has just the sources Mean and Lines, but the pseudofactors $L_1, \ldots, L_8$ refine this to the structure $\mathcal{R}_2$ whose sources are in the bottom diagram. Reduce $\mathcal{R}_2$ to $\mathcal{R}_3$ by letting $\text{Lines}_R$ be the sum of the sources for $L_1$, $L_3$, $L_5$ and $L_7$, and $\text{Lines}_T$ the sum of the sources for $L_2$, $L_4$, $L_6$ and $L_8$.

Now, $\mathcal{R}_1$ is structure balanced in relation to $\mathcal{Q}_1$, and $\mathcal{R}_3$ is structure balanced in relation to both $\mathcal{Q}_1$ and $\mathcal{Q}_2$. The three combined decompositions are in Table 4. In the third one, there are two sources involving $\text{Lines}_R$ because $\mathcal{R}_3$ is not orthogonal in relation to $\mathcal{Q}_2$. In fact, $\lambda_{P_1[B],\text{Lines}_R} = 1/4$. By Lemma 4.1, $P[B] \rhd \text{Lines}_R = \text{Lines}_R$. Because $P_1[B] \rhd \text{Lines}_R$ is a subspace of $P_1[B]$ and they have the same dimension, $P_1[B] \rhd \text{Lines}_R = P_1[B]$. However, $P_1[B] \neq \text{Lines}_R$, because the relevant efficiency factor is not 1. Lemma 4.2 shows that any structure on plots which is structure balanced in relation to $\mathcal{P}$ must be a refinement of $\mathcal{P}$. However, $P[B] \rhd \text{Lines}$ is not contained in any of the sources for $\mathcal{P}$; nor is $P[B] \rhd \text{Lines}_R$. Hence, of these three decompositions, only $\mathcal{Q}_2 \rhd \mathcal{R}_3$ is structure balanced in relation to $\mathcal{P}$. This is why the pseudofactors for lines and plots are chosen as they are.

The full decomposition is in Table 5.

The first two examples in [40] involve randomized-inclusive randomizations. This is not immediately obvious in [40], because the pseudofactors remain implicit, as do extra factors, such as Replicates, that refine $\mathcal{Q}$.

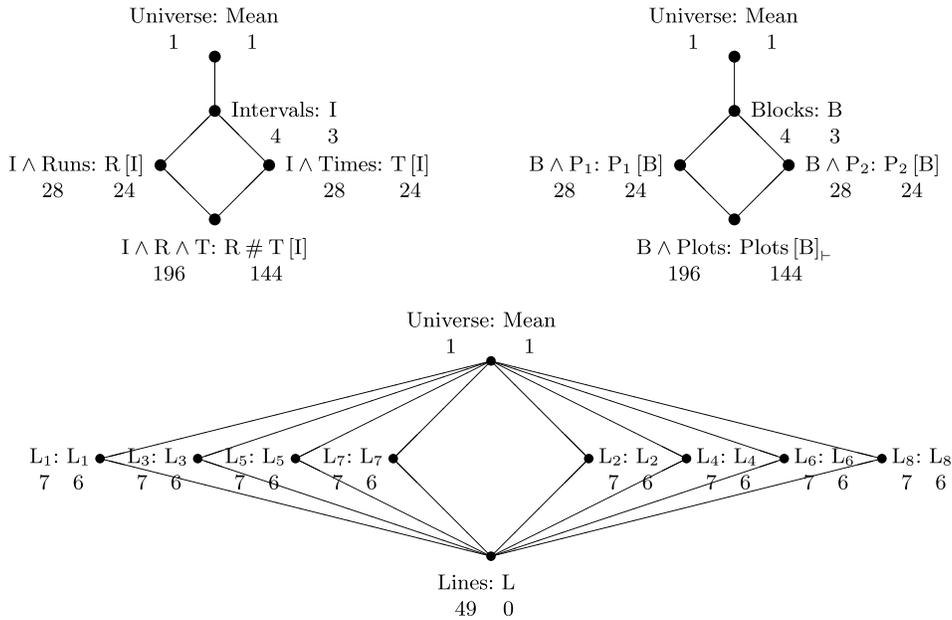

Fig. 12. *Hasse diagrams for analyses, plots and lines in Example 5.*

Table 4
*Three decompositions of the plots space in Example 5*

| $\mathcal{Q}_1 \rhd \mathcal{R}_1$ | | $\mathcal{Q}_1 \rhd \mathcal{R}_3$ | | $\mathcal{Q}_2 \rhd \mathcal{R}_3$ | |
|---|---|---|---|---|---|
| **source** | **d.f.** | **source** | **d.f.** | **source** | **d.f.** |
| Mean | 1 | Mean | 1 | Mean | 1 |
| Blocks | 3 | Blocks | 3 | Blocks | 3 |
| $P[B] \rhd \text{Lines}$ | 48 | $P[B] \rhd \text{Lines}_R$ | 24 | $P_1[B] \rhd \text{Lines}_R$ | 24 |
| $P[B] \vdash \mathcal{R}_1$ | 144 | $P[B] \rhd \text{Lines}_T$ | 24 | $P_2[B] \rhd \text{Lines}_T$ | 24 |
| | | $P[B] \vdash \mathcal{R}_3$ | 144 | $\text{Plots}[B]_\vdash \rhd \text{Lines}_R$ | 24 |
| | | | | $\text{Plots}[B]_\vdash \rhd \text{Lines}_T$ | 24 |
| | | | | $\text{Plots}[B]_\vdash \vdash \mathcal{R}_3$ | 96 |

Table 5
*Decomposition table for Example 5*

| analyses tier | | plots tier | | | lines tier | | |
|---|---|---|---|---|---|---|---|
| **source** | **d.f.** | **eff.** | **source** | **d.f.** | **eff.** | **source** | **d.f.** |
| Mean | 1 | 1 | Mean | 1 | 1 | Mean | 1 |
| Intervals | 3 | 1 | Blocks | 3 | | | |
| Runs[I] | 24 | 1 | $P_1[B]$ | 24 | $\frac{1}{4}$ | $\text{Lines}_R$ | 24 |
| Times[I] | 24 | 1 | $P_2[B]$ | 24 | $\frac{1}{4}$ | $\text{Lines}_T$ | 24 |
| R#T[I] | 144 | 1 | $\text{Plots}[B]_\vdash$ | 144 | $\frac{3}{4}$ | $\text{Lines}_R$ | 24 |
| | | | | | $\frac{3}{4}$ | $\text{Lines}_T$ | 24 |
| | | | | | | Residual | 96 |

**7. Structure-balanced experiments with a longer chain of randomizations.** It is clear that the results of Section 5 can be generalized by induction to experiments involving a chain of more than two randomizations so long as each successive pair is either composed or randomized-inclusive. For example, multiphase experiments with three or more phases will commonly employ three or more such randomizations. Suppose that $\Omega_i$ has structure $\mathcal{P}_i$, for $i = 1, \ldots, p$, and that the chain of randomizations goes from $\Omega_p$ to $\Omega_{p-1}$ to ... to $\Omega_1$. Assume that, if necessary, pseudofactors have been included in each structure, as in Section 6. Then the decomposition is

$$(\cdots (\mathcal{P}_1 \rhd \mathcal{P}_2) \rhd \cdots) \rhd \mathcal{P}_p \quad \text{or} \quad \mathcal{P}_1 \rhd (\mathcal{P}_2 \rhd (\cdots \rhd \mathcal{P}_p) \cdots).$$

Theorem 5.1 shows that the operator $\rhd$ is associative, and so the parentheses may be omitted.

The following example has four tiers and three randomizations. It appears quite simple, because, in each of the two proposed designs, the relationship between each successive pair of structures is orthogonality. Nevertheless, it shows how examination of the decomposition tables can aid the experimenter to choose among competing designs.

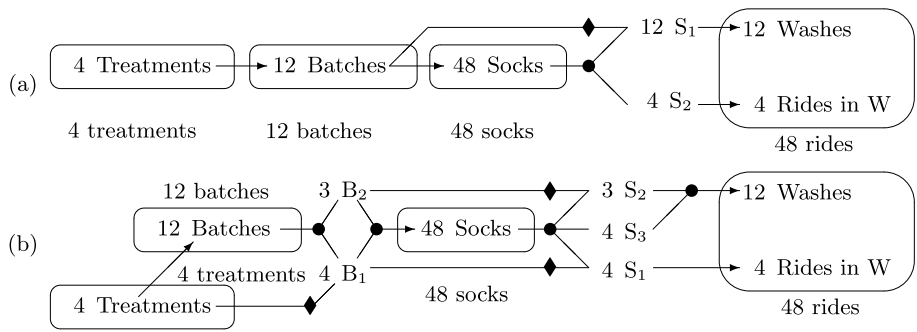

Fig. 13. *Two methods of randomization in Example 6: socks are randomly split into batches, treatments are randomized to batches, and socks are randomized to rides;* $S_1, \ldots, S_3$ *are pseudofactors for Socks (different in the two methods);* $B_1$ *and* $B_2$ *are pseudofactors for Batches.*

EXAMPLE 6 (Knitted socks). Example 14 in [14] concerns an experiment in which 48 knitted socks are randomly divided into twelve batches, to which four chlorination treatments are randomized. This develops Example 5.9 from [17], in which various ways of randomizing the socks to rides in washing machines are suggested. Figure 13 shows two of these. In both cases, the first and second randomizations are composed while the second and third are randomized-inclusive. The results of Sections 5 and 6 can be applied directly to give the decompositions in Tables 6 and 7. In Table 6, the source Socks $\vdash S_1$ is the part of Socks that is orthogonal to the source $S_1$. The source Socks $\vdash (S_2 \wedge S_3)$ in Table 7 is defined similarly. The first experiment has the advantage that the Residual for Treatments has more degrees of freedom but the second has the advantage that the Treatments source is confounded with Rides[Washes] rather than the potentially more variable Washes.

**8. Summary.** In this paper, the notation $\mathcal{P} \triangleright \mathcal{Q}$ has been introduced for the orthogonal decomposition of $V_\Omega$ for a single randomization, based on the two structures $\mathcal{P}$ and $\mathcal{Q}$. The idea has been extended to cover composed and randomized-inclusive randomizations, which involve three structures, $\mathcal{P}$, $\mathcal{Q}$ and $\mathcal{R}$. For these two related types of multiple randomization, it has been proven that the orthogonal decompositions $(\mathcal{P} \triangleright \mathcal{Q}) \triangleright \mathcal{R}$ and $\mathcal{P} \triangleright (\mathcal{Q} \triangleright \mathcal{R})$ of $V_\Omega$ are the same. As a result, the decomposition may be done from left-to-right or right-to-left. The decomposition table for a left-to-right decomposition, when there are three tiers, can be produced using the `AMTIER` procedure in GenStat [16].

TABLE 6
*Decomposition table for Example 6: first method*

| rides tier | | socks tier | | batches tier | | treatments tier | |
|---|---|---|---|---|---|---|---|
| **source** | **d.f.** | **source** | **d.f.** | **source** | **d.f.** | **source** | **d.f.** |
| Mean | 1 | Mean | 1 | Mean | 1 | Mean | 1 |
| Washes | 11 | $S_1$ | 11 | Batches | 11 | Treatments | 3 |
| | | | | | | Residual | 8 |
| Rides[Washes] | 36 | Socks $\vdash S_1$ | 36 | | | | |

TABLE 7
*Decomposition table for Example 6: second method*

| rides tier | | socks tier | | batches tier | | treatments tier | |
|---|---|---|---|---|---|---|---|
| **source** | **d.f.** | **source** | **d.f.** | **source** | **d.f.** | **source** | **d.f.** |
| Mean | 1 | Mean | 1 | Mean | 1 | Mean | 1 |
| Washes | 11 | $S_2 \wedge S_3$ | 11 | $B_2$ | 2 | | |
| | | | | Residual | 9 | | |
| Rides[Washes] | 36 | Socks $\vdash (S_2 \wedge S_3)$ | 36 | Batches $\vdash B_2$ | 9 | Treatments | 3 |
| | | | | | | Residual | 6 |
| | | | | Residual | 27 | | |

When it comes to designing experiments and evaluating these designs, we have obtained the particularly useful result that, if each design is structure balanced, then so is the composite. Further, its matrix of efficiency factors is the product of those for the individual designs.

The question remains as to how the results for multiple randomizations that are unrandomized-inclusive, independent, coincident or double [14] compare with those presented here. It will be addressed in a further paper.

**Acknowledgments.** The authors are grateful to the referees for suggestions that helped to improve the readability of the paper.

SCHOOL OF MATHEMATICS AND STATISTICS  
UNIVERSITY OF SOUTH AUSTRALIA  
NORTH TERRACE, ADELAIDE, SA 5000  
AUSTRALIA  
E-MAIL: chris.brien@unisa.edu.au  
URL: http://chris.brien.name  

SCHOOL OF MATHEMATICAL SCIENCES  
QUEEN MARY, UNIVERSITY OF LONDON  
MILE END ROAD  
LONDON E1 4NS  
UNITED KINGDOM  
E-MAIL: r.a.bailey@qmul.ac.uk